\documentclass[12pt]{article}

\usepackage{amssymb,bm}
\usepackage{amsmath}
\usepackage{graphicx}

\newtheorem{theorem}{{\sc Theorem}}

\newtheorem{lemma}{{\sc Lemma}}
\newtheorem{cor}{Corollary}

\newtheorem{defn}{Definition}

\newcommand{\rf}[1]{(\ref{#1})}

\newcommand{\bea}[1]{\begin{array}{#1}}
\newcommand{\ea}{\end{array}}

\newcommand{\Var}{\mathrm{Var}}

\newcommand{\wk}[1]{{\buildrel #1\over\rightharpoonup}\:}
\newcommand{\weak}{\rightharpoonup\:}
\newcommand{\dOm}{\partial\Omega}

\newcommand{\bra}[1]{\overline{#1}}

\newcommand{\Trc}{\mathrm{Tr}\,}

\newcommand{\hf}{\displaystyle\frac{1}{2}}
\newcommand{\nth}[1]{\displaystyle\frac{1}{#1}}

\newcommand{\Grad}{\nabla}
\newcommand{\Div}{\nabla \cdot}

\newcommand{\Md}{\partial}

\renewcommand{\Hat}[1]{\widehat{#1}}
\newcommand{\Tld}[1]{\widetilde{#1}}

\newcommand{\av}[1]{\langle #1 \rangle}

\newenvironment{proof}[1][Proof]{\trivlist\item[\hskip\labelsep{\hskip    
\parindent{\normalfont\scshape#1:}\hskip .321429\parindent}]\ignorespaces}
{\qed\endtrivlist}
\newcommand\qed{\ \rule[-0.2ex]{0.3em}{1.5ex}}

\newcommand{\limi}{\mathop{\underline\lim}}

\newcommand{\bc}{boundary condition}
\newcommand{\rhs}{right hand side}

\newcommand{\nbh}{neighborhood}
\newcommand{\IFF}{if and only if }

\newcommand{\beqa}{\begin{eqnarray}}
\newcommand{\beqn}{\begin{eqnarray*}}
\newcommand{\eeqa}[1]{\label{#1}\end{eqnarray}}
\newcommand{\eeqn}{\end{eqnarray*}}
\newcommand{\beq}{\begin{equation}}
\newcommand{\eeq}[1]{\label{#1}\end{equation}}

\newcommand{\bb}[1]{\mathbb{ #1}}

\newcommand{\Ga}{\alpha}
\newcommand{\Gb}{\beta}
\newcommand{\Gd}{\delta}
\newcommand{\Ge}{\epsilon}

\newcommand{\Gc}{\chi}

\newcommand{\Gl}{\lambda}

\newcommand{\Gm}{\mu}
\newcommand{\Gv}{\nu}
\newcommand{\Gp}{\pi}

\newcommand{\Gth}{\theta}

\newcommand{\Gu}{\upsilon}

\newcommand{\Gx}{\xi}

\newcommand{\GD}{\Delta}

\newcommand{\GL}{\Lambda}

\newcommand{\GO}{\Omega}

\bmdefine\BGa{\alpha}
\bmdefine\BGb{\beta}
\bmdefine\BGd{\delta}
\bmdefine\BGe{\epsilon}
\bmdefine\BGve{\varepsilon}
\bmdefine\BGf{\phi}
\bmdefine\BGvf{\varphi}
\bmdefine\BGg{\gamma}
\bmdefine\BGc{\chi}
\bmdefine\BGi{\iota}
\bmdefine\BGk{\kappa}
\bmdefine\BGl{\lambda}
\bmdefine\BGn{\eta}
\bmdefine\BGm{\mu}
\bmdefine\BGv{\nu}
\bmdefine\BGp{\pi}
\bmdefine\BGth{\theta}
\bmdefine\BGvth{\vartheta}
\bmdefine\BGr{\rho}
\bmdefine\BGvr{\varrho}
\bmdefine\BGs{\sigma}
\bmdefine\BGvs{\varsigma}
\bmdefine\BGt{\tau}
\bmdefine\BGj{\tau}
\bmdefine\BGu{\upsilon}
\bmdefine\BGo{\omega}
\bmdefine\BGx{\xi}
\bmdefine\BGy{\psi}
\bmdefine\BGz{\zeta}
\bmdefine\BGD{\Delta}
\bmdefine\BGF{\Phi}
\bmdefine\BGG{\Gamma}
\bmdefine\BGL{\Lambda}
\bmdefine\BGP{\Pi}
\bmdefine\BGT{\Theta}
\bmdefine\BGS{\Sigma}
\bmdefine\BGU{\Upsilon}
\bmdefine\BGO{\Omega}
\bmdefine\BGX{\Xi}
\bmdefine\BGY{\Psi}

\newcommand{\CA}{{\cal A}}
\newcommand{\CB}{{\cal B}}

\newcommand{\CF}{{\cal F}}

\newcommand{\CI}{{\cal I}}

\newcommand{\CM}{{\cal M}}

\newcommand{\CO}{{\cal O}}

\newcommand{\CR}{{\cal R}}
\newcommand{\CS}{{\cal S}}

\newcommand{\CY}{{\cal Y}}

\bmdefine\BCa{{\cal a}}
\bmdefine\BCb{{\cal b}}
\bmdefine\BCc{{\cal c}}
\bmdefine\BCd{{\cal d}}
\bmdefine\BCe{{\cal e}}
\bmdefine\BCf{{\cal f}}
\bmdefine\BCg{{\cal g}}
\bmdefine\BCh{{\cal h}}
\bmdefine\BCi{{\cal i}}
\bmdefine\BCj{{\cal j}}
\bmdefine\BCk{{\cal k}}
\bmdefine\BCl{{\cal l}}
\bmdefine\BCm{{\cal m}}
\bmdefine\BCn{{\cal n}}
\bmdefine\BCo{{\cal o}}
\bmdefine\BCp{{\cal p}}
\bmdefine\BCq{{\cal q}}
\bmdefine\BCr{{\cal r}}
\bmdefine\BCs{{\cal s}}
\bmdefine\BCt{{\cal t}}
\bmdefine\BCu{{\cal u}}
\bmdefine\BCv{{\cal v}}
\bmdefine\BCx{{\cal x}}
\bmdefine\BCy{{\cal y}}
\bmdefine\BCz{{\cal z}}
\bmdefine\BCA{{\cal A}}
\bmdefine\BCB{{\cal B}}
\bmdefine\BCC{{\cal C}}
\bmdefine\BCD{{\cal D}}
\bmdefine\BCE{{\cal E}}
\bmdefine\BCF{{\cal F}}
\bmdefine\BCG{{\cal G}}
\bmdefine\BCH{{\cal H}}
\bmdefine\BCI{{\cal I}}
\bmdefine\BCJ{{\cal J}}
\bmdefine\BCK{{\cal K}}
\bmdefine\BCL{{\cal L}}
\bmdefine\BCM{{\cal M}}
\bmdefine\BCN{{\cal N}}
\bmdefine\BCO{{\cal O}}
\bmdefine\BCP{{\cal P}}
\bmdefine\BCQ{{\cal Q}}
\bmdefine\BCR{{\cal R}}
\bmdefine\BCS{{\cal S}}
\bmdefine\BCT{{\cal T}}
\bmdefine\BCU{{\cal U}}
\bmdefine\BCV{{\cal V}}
\bmdefine\BCW{{\cal W}}
\bmdefine\BCX{{\cal X}}
\bmdefine\BCY{{\cal Y}}
\bmdefine\BCZ{{\cal Z}}

\bmdefine\Bzr{ 0}
\bmdefine\Ba{ a}
\bmdefine\Bb{ b}
\bmdefine\Bc{ c}
\bmdefine\Bd{ d}
\bmdefine\Be{ e}
\bmdefine\Bf{ f}
\bmdefine\Bg{ g}
\bmdefine\Bh{ h}
\bmdefine\Bi{ i}
\bmdefine\Bj{ j}
\bmdefine\Bk{ k}
\bmdefine\Bl{ l}
\bmdefine\Bm{ m}
\bmdefine\Bn{ n}
\bmdefine\Bo{ o}
\bmdefine\Bp{ p}
\bmdefine\Bq{ q}
\bmdefine\Br{ r}
\bmdefine\Bs{ s}
\bmdefine\Bt{ t}
\bmdefine\Bu{ u}
\bmdefine\Bv{ v}
\bmdefine\Bw{ w}
\bmdefine\Bx{ x}
\bmdefine\By{ y}
\bmdefine\Bz{ z}
\bmdefine\BA{ A}
\bmdefine\BB{ B}
\bmdefine\BC{ C}
\bmdefine\BD{ D}
\bmdefine\BE{ E}
\bmdefine\BF{ F}
\bmdefine\BG{ G}
\bmdefine\BH{ H}
\bmdefine\BI{ I}
\bmdefine\BJ{ J}
\bmdefine\BK{ K}
\bmdefine\BL{ L}
\bmdefine\BM{ M}
\bmdefine\BN{ N}
\bmdefine\BO{ O}
\bmdefine\BP{ P}
\bmdefine\BQ{ Q}
\bmdefine\BR{ R}
\bmdefine\BS{ S}
\bmdefine\BT{ T}
\bmdefine\BU{ U}
\bmdefine\BV{ V}
\bmdefine\BW{ W}
\bmdefine\BX{ X}
\bmdefine\BY{ Y}
\bmdefine\BZ{ Z}

\newcommand{\SFI}{\mbox{\sf I}}

\newcommand{\SFL}{\mbox{\sf L}}

\title{Direct approach to the problem of strong local minima 
  in Calculus of Variations}
\author{Yury Grabovsky \and Tadele Mengesha}

\begin{document}

\maketitle
\begin{abstract}
  The paper introduces a general strategy for identifying strong local
  minimizers of variational functionals.  It is based on the idea that any
  variation of the integral functional can be evaluated directly in terms of
  the appropriate parameterized measures.  We demonstrate our approach on a
  problem of $W^{1,\infty}$ weak-* local minima---a slight weakening of the
  classical notion of strong local minima. We obtain the first
  quasiconvexity-based set of sufficient conditions for $W^{1,\infty}$ weak-*
  local minima.
\end{abstract}

\section{Introduction}
\setcounter{equation}{0}
\label{section:intro}
\label{sub:abh}
In this paper we consider the class of integral functionals of the form
\beq
E(\By)=\int_{\GO}W(\Bx,\Grad\By(\Bx))d\Bx,
\eeq{energy}
where $\GO$ is a smooth (i.e. of class $C^{1}$) and bounded domain in $\bb{R}^{d}$
and the Lagrangian $W:\bra{\GO}\times\bb{R}^{m\times d}\to\bb{R}$ is assumed to
be a continuous function. The symbol $\bb{R}^{m\times d}$ is used to denote
the space of all $m\times d$ real matrices. The functional \rf{energy} is
defined on the set of admissible functions
\beq
\CA=\{\By\in W^{1,\infty}(\GO;\bb{R}^{m}): 
\By(\Bx)=\Bg(\Bx),\ \Bx\in\bra{\Md\GO_{1}}\},
\eeq{bc0}
where $\dOm_{1}$ and $\dOm_{2}=\dOm\setminus\bra{\dOm_{1}}$ are smooth 
(i.e. of class $C^{1}$) relatively open subsets of $\dOm$, and 
$\Bg\in C^{1}(\bra{\dOm_{1}};\bb{R}^{m})$. We omit the dependence of $W$ on
$\By$ to simplify our analysis and because such dependence does not introduce
conceptually new difficulties 
(within the context of our discussion). The omission of dependence of $W$ on $\Bx$,
however, does not lead to similar simplifications, as the dependence on $\Bx$
will reappear in our analysis even if $W$ does not depend on $\Bx$ explicitly.

A fundamental problem in Calculus of Variations and its applications is the
problem of finding local minimizers (see \cite[Problem 9]{ball02}, for example). 
The notion of the local minimizer, in contrast to the global one, depends in an essential 
way on the topology on the space $\CA$ of functions on which the variational
functional is defined. We assume that the topology on $\CA$ comes from a
topological vector space topology $\tau$ on $W^{1,\infty}(\GO;\bb{R}^{m})$, since we
want standard linear operations to be continuous.
Let
\beq
\Var(\CA)=\{\BGf\in W^{1,\infty}(\GO;\bb{R}^{m}): \BGf|_{\Md\GO_{1}} =\Bzr\}
\eeq{vars}
be the space of variations. Observe that $\By+\BGf\in\CA$ for all $\By\in\CA$ and all
$\BGf\in\Var(\CA)$.
\begin{defn}
  \label{def:tauvar}
The sequence $\{\BGf_{n}: n\ge 1\}\subset\Var(\CA)$ is called a
\textbf{$\tau$-variation} if $\BGf_{n}\to\Bzr$ in $\tau$.
\end{defn}
\begin{defn}
\label{def:locmin}
We say that $\By\in\CA$ is a \textbf{$\tau$-local minimum}, if for every
$\tau$-variation $\{\BGf_{n}: n\ge 1\}\subset\Var(\CA)$ there exists $N\ge 1$ 
such that $E(\By)\le E(\By+\BGf_{n})$ for all $n\ge N$.
\end{defn}
The classical notions of strong and weak local minima are examples of
$\tau$-local minima, where $\tau$ is the $L^{\infty}$ and $W^{1,\infty}$
topologies on $W^{1,\infty}(\GO;\bb{R}^{m})$ respectively. Clearly, the weaker the
topology $\tau$, the stronger the notion of the local minimum. This is
reflected in the terminology. The notion of strong local minimum is
stronger than the notion of the weak one. 
\begin{defn}
  \label{def:wkstr}
A variation is called \textbf{strong} or \textbf{weak} if it is an
$L^{\infty}$ variation or a $W^{1,\infty}$ variation respectively.
\end{defn}

If the topology $\tau$ is non-metrizable, like the $W^{1,\infty}$ weak-*
topology considered in this paper, then the sequence-based definition is
different from the one based on open sets. In this paper we will use
the sequence-based Definition~\ref{def:locmin}. 

The problem of strong local minima is fairly well-understood in the classical
Calculus of Variations, $d=1$ (Weierstrass) or $m=1$ (Hestenes \cite{hest48}).
The present paper will focus on the case $d>1$ and $m>1$, where many
fundamental problems still remain open largely because the existing methods
are not as effective in this case as they are in the classical cases. In this
paper we bring the analytical machinery developed for the ``Direct Method'' in
Calculus of Variations, introduced by Tonelli for the purpose of proving
existence of \emph{global} minimizers, to bear on the problem of \emph{local}
minimizers. We propose a general strategy that is capable of delivering
quasiconvexity-based sufficient conditions for strong local minima. We
demonstrate how our strategy works in a simplified setting of smooth (i.e.
$C^{1}$) extremals $\By(\Bx)$ and stronger (i.e. $W^{1,\infty}$ weak-*)
topology $\tau$.  Strengthening topology $\tau$
from $L^{\infty}$ to $W^{1,\infty}$ weak-* means that we restrict possible
variations $\{\BGf_{n}\}$ to sequences that converge to zero uniformly, while
remaining bounded in $W^{1,\infty}(\GO;\bb{R}^{m})$. In other words, the
$W^{1,\infty}$ weak-* variations are the sequences that converge to zero
$W^{1,\infty}$ weak-*. From this point on the word ``variation'' will mean
$W^{1,\infty}$ weak-* variation.

Our approach should also be
applicable even if $\By(\Bx)$ is not of class  
$C^{1}$ and the topology $\tau$ is $L^{\infty}$. However, the actual technical
implementation will require overcoming a set of difficulties related to
the appearance of new necessary conditions on 
the behavior of $W$ at the discontinuities of
$\Grad\By(\Bx)$ and at infinity (see \cite{grtr08} for details). 

So far we did not require that the Lagrangian $W(\Bx,\BF)$ be
smooth. We do not want to make a global smoothness assumption on $W$ in order not to
rule out examples where the Lagrangian is piecewise smooth. For example, in
the mathematical theory of composite materials or optimal design the
Lagrangian is given as a minimum of finitely many quadratic functions
\cite{ks}. In fact, we do not need the Lagrangian $W$ to be smooth everywhere.
Let
\[
\CR = \{\BF \in \bb{R}^{m\times d}: \BF = \Grad\By (\Bx) ~\textrm{for some}~ \Bx\in
\bra{\GO}\}.
\]
In other words, $\CR$ is the range of
$\Grad\By(\Bx)$.  We assume that $W$ is of class $C^{2}$ on $\CR$, meaning that
there exists an open set $\CO$ such that $\CR\subset \CO$ and the functions
$W(\Bx,\BF),$ $W_{\BF}(\Bx,\BF),$ and $W_{\BF\BF}(\Bx,\BF)$ are continuous on
$\bra{\GO}\times\CO$. Throughout the paper we will use the subscript notation 
to denote the vectors, matrices and higher order tensors of partial derivatives.


\section{The strategy for identifying strong local minima}
\setcounter{equation}{0}
\label{sec:fp}
One of the fundamental problems of Calculus of Variations is to find
sufficient conditions for strong local minima.  This problem (for $d>1$ and
$m>1$) is quite old and there are many sets of sufficient conditions that have
already been found \cite{cara29,deDon35,lepa41,tahe01,weyl35}. However, none
of them is in any sense close to the necessary conditions that are formulated
using the notion of quasiconvexity. In recent years
it became clear, that the the quasiconvexity condition is the correct
multi-dimensional analog of the classical Weierstrass condition (positivity of
the Weierstrass excess function) \cite{bama84}.  The quasiconvexity
condition was first introduced by Morrey \cite{morr52}, who showed that this
condition is necessary and sufficient for $W^{1,\infty}$ weak-* lower
semicontinuity of the variational integrals \rf{energy}. 

In this paper we
present the first set of \emph{quasiconvexity based} sufficient conditions for
$W^{1,\infty}$ weak-* local minima.
Our strategy is the result of the insights achieved in \cite{grtr08},
where the necessary conditions for strong local minima are examined in greater
generality. In this paper we will only need the observation made in
\cite{grtr08} that the limit 
\beq
\Gd E=\limi_{n\to\infty}\frac{\GD E(\BGf_{n})}{\|\Grad\BGf_{n}\|_{2}^{2}},
\eeq{rightstuff}
where
\[
\GD E(\BGf_{n})=\int_{\GO}(W(\Grad\By+\Grad\BGf_{n})-W(\Grad\By))d\Bx,
\]
is always finite for an extremal $\By(\Bx)$ (i.e. solution of the Euler-Lagrange
equation \rf{e-l eqn} below).
Moreover, the requirement of non-negativity of $\Gd E$ for 
specific variations $\BGf_{n}$ produces all known necessary conditions for a $C^{1}$
extremal $\By(\Bx)$ to be a strong local minimizer.  
In \rf{rightstuff} and throughout the paper $\|\BA\|_{p}$ denotes the $L^{p}$ norm of
the Frobenius norm $|\BA(\Bx)|$ of the matrix field $\BA(\Bx)$.

We remark that the choice of the denominator in \rf{rightstuff} is not
arbitrary. It expresses the correct size scale of the increment of the functional 
under the variation $\BGf_{n}$. 
Now, we are ready to describe our strategy for identifying strong local minima.
\begin{center}
  \textbf{The strategy:}
\end{center}
\label{strat}
\begin{enumerate}
\item[\textbf{Step 1.}] Make specific variations for which $\Gd E$ can be
  computed explicitly.  Obtain
necessary conditions for $\By\in\CA$ to be a strong local minimizer from the
  inequality $\Gd E\ge 0$. 
\item[\textbf{Step 2.}] Prove that if $\By(\Bx)$ satisfies
the necessary conditions from Step 1, then $\Gd E\ge 0$ for \emph{all} variations
$\{\BGf_{n}\}$.  
\item[\textbf{Step 3.}] Characterize those variations $\{\BGf_{n}\}$ for which
$\Gd E=0$. 
\item[\textbf{Step 4.}] Formulate the weakest additional conditions, that
  together with the necessary conditions obtained on Step 1, prevent $\GD
  E(\BGf_{n})$ from becoming negative for large $n$ for variations,
  for which $\Gd E=0$. 
\end{enumerate}
In Step 1, the necessary conditions for $C^{1}$ functions $\By(\Bx)$ are
well-known by now. They consist of the Euler-Lagrange equation, non-negativity of
second variation and the quasiconvexity conditions in the interior and on the
free boundary \cite{bama84}. For more general Lipschitz extremals $\By(\Bx)$
other necessary conditions may appear (see \cite{grtr08} for a discussion of
why this happens). 
Step 2 is the focus of the present paper.  Step 4 should naturally follow from
the analysis of Step 3. At this moment Step 3 is still open.  We avoid the
delicate analysis entailed by Step 3 by imposing extra 
conditions that prevent \emph{any} non-zero variation to satisfy $\Gd E=0$.


\section{Reformulation of the problem}
\setcounter{equation}{0}
\label{sec:ref}
Our first observation is that the Euler-Lagrange equation 
\begin{equation}
\label{e-l eqn}
\left\{ \begin{array}{ll}
\Div W_{\BF}(\Bx,\BF (\Bx)) = \Bzr, & \Bx\in\GO,\\[1ex]
W_{\BF}(\Bx,\BF (\Bx))\Bn (\Bx) = \Bzr, & \Bx\in\dOm _{2},
\end{array} \right.
\end{equation}
where $\Bn(\Bx)$ is the outer unit normal to $\dOm$ at $\Bx\in\dOm$,
can be 
completely decoupled from the other necessary conditions for strong local minima.
This is done by replacing the functional increment $\GD E(\BGf_{n})$ by
\beq
\GD' E (\BGf_{n})=  \int_{\GO}W^{\circ}(\Bx,\Grad\BGf_{n}(\Bx))d\Bx,
\eeq{Dpr}
where
\beq
W^{\circ}(\Bx,\BF)=
W(\Bx,\BF (\Bx)+\BF) - W(\Bx,\BF(\Bx)) - (W_{\BF}(\Bx,\BF(\Bx)),\BF)
\eeq{Wcirc}
is related to the Weierstrass excess function.
In the formula above and throughout the paper we use the notation
$\BF(\Bx)=\Grad\By(\Bx)$
and the inner product notation $(\cdot,\cdot)$ corresponding to the dot
product on $\bb{R}^{d}$ and the Frobenius inner product
$(\BA,\BB)=\Trc(\BA\BB^{T})$ on $\bb{R}^{m\times d}.$ 

We conclude, therefore, that the role of
the Euler-Lagrange equation \rf{e-l eqn} is to establish equivalence between 
$\GD'E(\BGf_{n})$---a quantity that our analysis applies to, and the functional 
increment $\GD E(\BGf_{n})$---a quantity with variational meaning.
We can view the transition from $\GD E$ to $\GD'E$ as a transformation
\beq
\Pi:(W(\Bx,\BF),\By(\Bx))\mapsto (W^{\circ}(\Bx,\BF),\Bzr).
\eeq{transf} 

We note, that regardless of the choice of $\By(\Bx)$, the function
$\Bzr$ satisfies the Euler-Lagrange equation for the Lagrangian
$W^{\circ}$. Moreover, it is clear, that $\By(\Bx)$ is a $\tau$-local minimum
for the Lagrangian $W$ \IFF $\By(\Bx)$ solves the Euler-Lagrange equation
\rf{e-l eqn} and $\Bzr$ is a $\tau$-local minimum for the functional
with Lagrangian $W^{\circ}$, since the functional increment $\GD E$ for 
$W^{\circ}$ is exactly
$\GD' E$ for $W$.  Thus, the projection $\Pi$ given by \rf{transf}, (it
is easy to verify that $\Pi$ is indeed a projection) allows us to decouple
the Euler-Lagrange equation from all the other conditions that one would
require to guarantee that $\By(\Bx)$ is a local minimizer. The range of $\Pi$
is a set of continuous
functions $W^{\circ}(\Bx,\BF)$ that are twice continously differentiable on
some \nbh\ of
$\BF=\Bzr$ and vanish with its first derivative at $\BF=\Bzr$.  It will be
convenient for us to represent $W^{\circ}$ in the form that shows the
quadratic term in its Taylor expansion around $\BF=\Bzr$ explicitly, because
it appears  in the formula for the second variation. 
\beq
W^{\circ}(\Bx,\BF)=\hf(\SFL(\Bx)\BF,\BF)+|\BF|^{2}U(\Bx,\BF),
\eeq{UL}
where
\beq
\SFL(\Bx)=W^{\circ}_{\BF\BF}(\Bx,\Bzr)=W_{\BF\BF}(\Bx,\BF(\Bx))
\eeq{Ldef}
and
\[
U(\Bx,\BF)=\nth{|\BF|^{2}}\left(W^{\circ}(\Bx,\BF) - 
\hf(\SFL(\Bx)\BF,\BF)\right)
\]
is a continuous function on $\bra{\GO} \times\bb{R}^{m\times d}$ that vanishes
on $\bra{\GO} \times\{\Bzr\}$.

Replacing $\GD E$ with $\GD' E$ and $W$ with $W^{\circ}$, we reduce the
problem of local minima
to the determination of the sign of $\Gd' E$ given by 
\beq
\Gd' E=\limi_{n\to\infty}\frac{\GD' E(\BGf_{n})}{\|\Grad\BGf_{n}\|_{2}^{2}}
=\limi_{n\to\infty}\nth{\|\Grad\BGf_{n}\|_{2}^{2}}\int_{\GO}
W^{\circ}(\Bx,\Grad\BGf_{n})d\Bx.
\eeq{dpr}
We reiterate that  $\Gd'E=\Gd E$ for all variations $\BGf_{n}$ \IFF $\By(\Bx)$ 
satisfies the Euler-Lagrange equation \rf{e-l eqn}. Substituting  the
representation \rf{UL} of $W^{\circ}$ into \rf{dpr}, we obtain
\beq 
\Gd'E=\limi_{n\to\infty}
\int_{\GO}\left(U(\Bx, \Ga_{n}\Grad \BGy_{n}(\Bx))|\Grad \BGy _{n}(\Bx)|^{2} +
\hf(\SFL(\Bx)\Grad \BGy_{n}(\Bx), \Grad \BGy_{n}(\Bx))\right )d\Bx, 
\eeq{increment}
where 
\beq
\Ga_{n} = \|\Grad \BGf_{n}\|_{2}\text{ and }
\BGy_{n}(\Bx)=\frac{\BGf_{n}(\Bx)}{\|\Grad\BGf_{n}\|_{2}}.
\eeq{alpsi}
The formula \rf{increment} will serve as a starting point of our analysis. 
In order to simplify notation we will use a shorthand
\beq
\CF (\Bx, \Ga, \BG)=\frac{W^{\circ}(\Bx,\Ga\BG)}{\Ga^{2}}
=U(\Bx,\Ga\BG)|\BG|^{2} + \hf(\SFL(\Bx)\BG, \BG).
\eeq{shhd}
Thus, in terms of $\CF$
\beq
\Gd' E =
\limi_{n\to \infty}\int_{\GO}\CF (\Bx, \Ga_{n}, \Grad\BGy_{n})d\Bx.
\eeq{limit-divided-increment}

Finally, we would like to note that our approach is in some sense dual to the
classical approach that studies the effect of a family of variations on a
given integral functional. Borrowing the idea of duality from Young
\cite{young33,young37} (see also the papers \cite{ball89,tart79a} that helped
bring the importance of Young measures for applications), 
we consider a given variation $\BGf_{n}$ and study its
effect on pairs $(U(\Bx,\BF),Q(\Bx,\BF))$, where $U$ varies in the space of
continuous functions on $\bra{\GO}\times\bb{R}^{m\times d}$ that vanish at
$\bra{\GO}\times\{\Bzr\}$ and $Q(\Bx,\BF)=(\SFL(\Bx)\BF, \BF)$ is quadratic in
$\BF$ and continuous in $\Bx\in\bra{\GO}$.

The formula \rf{increment} indicates that we prefer to regard a variation
$\{\BGf_{n}\}$ as a pair $(\Ga_{n},\BGy_{n})$, where $\Grad\BGy_{n}$ has
$L^{2}$-norm equal to 1 and
$\Ga_{n}\Grad\BGy_{n}(\Bx)$ is bounded in $L^{\infty}$. We can think of
$\Ga_{n}$ as the ``size'' of the variation and of $\BGy_{n}$ as its ``shape''.


\section{Necessary conditions and sufficient conditions}
\label{sub:nec}
\setcounter{equation}{0}
We begin with a quick recap of the known necessary conditions for strong local minima
for $\By\in C^{1}(\bra{\GO};\bb{R}^{m})$ (see, for
example, \cite{bama84}). We then show that necessary
conditions imply non-negativity of $\Gd' E$. Finally, we show that if we
strengthen the non-strict inequalities appearing in the necessary
conditions below, we will obtain sufficient conditions for $W^{1,\infty}$ weak-*
local minimizers of class $C^{1}$. (See Theorem~\ref{semisuff2} below.)

It is well-known that if we perturb $\By(\Bx)$ using
special weak variations 
\beq
\By(\Bx)\to\By(\Bx)+\Ge\BGf(\Bx),
\eeq{wkvar} 
we obtain the Euler-Lagrange equation \rf{e-l eqn} \emph{and} the condition of
non-negativity of the second variation 
\beq
\Gd^{2}E=\int_{\GO}(\SFL(\Bx)\Grad\BGf(\Bx),\Grad\BGf(\Bx))d\Bx
\eeq{secvar}
for all $\BGf\in\Var(\CA)$, where $\Var(\CA)$ is given by \rf{vars} and
$\SFL(\Bx)$ is given by \rf{Ldef}.

If we perturb $\By(\Bx)$ using the generalized ``Weierstrass needle''
\beq
\By(\Bx)\to\By(\Bx)+\Ge\BGf\left(\frac{\Bx-\Bx_{0}}{\Ge}\right),
\eeq{nuc}
where $\BGf(\Bx)\in W_{0}^{1,\infty}(B(\Bzr,1);\bb{R}^{m})$,
we will get the two quasiconvexity conditions: the Morrey
quasiconvexity condition \cite{morr52}
\beq
\int_{B(\Bzr,1)}W(\Bx_{0},\BF(\Bx_{0})  + \Grad \BGf (\Bx))d\Bx
\geq  \int_{B(\Bzr,1)} W(\Bx_{0},\BF (\Bx_{0}))d\Bx,
\eeq{qxi}
for all $\Bx_{0}\in\GO$, and the quasiconvexity at the free boundary condition
\cite{bama84}
\beq
\int _{B^{-}_{\Bn(\Bx_{0})}(\Bzr,1)}W(\Bx_{0},\BF(\Bx_{0})  + \Grad \BGf (\Bx))d\Bx \geq  
\int_{B^{-}_{\Bn(\Bx_{0})}(\Bzr,1)} W(\Bx_{0},\BF(\Bx_{0}))d\Bx,
\eeq{qxb}
for all $\Bx_{0}\in\dOm_{2}$. Here $B(\Bx,r)$ denotes an open ball in $\bb{R}^{d}$
centered at $\Bx$ with radius $r$ and $B^{-}_{\Bn}(\Bzr,1)$ denotes the half-ball
$
B^{-}_{\Bn}(\Bzr,1)=\{\Bx\in B(\Bzr,1),\ (\Bx,\Bn)<0\},
$ 
whose outer unit normal at the ``flat'' part of its boundary is equal to
$\Bn$. 

Morrey himself derived the quasiconvexity condition \rf{qxi} as a
necessary and sufficient condition for $W^{1,\infty}$ weak-* lower
semicontinuity of the integral functionals \rf{energy}. The necessity of
\rf{qxi} for strong local minimizers via the variation \rf{nuc} is due to Ball
and Marsden \cite{bama84}, even though the fact itself can be inferred from
the arguments of Meyers \cite{meye65}, whose focus was on lower semicontinuity
of integral functionals involving higher derivatives of $\By$. In fact, the
proof of Meyers' Lemma 1 in \cite{meye65} can be interpreted as a direct link
between $W^{1,\infty}$ weak-* local minima and $W^{1,\infty}$ weak-* lower
semicontinuity, explaining why Morrey's quasiconvexity appears naturally in both
contexts.

Our idea was to replace the original Lagrangian $W$ with the ``reduced Lagrangian''
$W^{\circ}(\Bx,\BF)$, given by \rf{Wcirc}. Therefore, we 
rewrite the quasiconvexity conditions \rf{qxi}--\rf{qxb} in terms of the
``reduced Lagrangian'' $W^{\circ}(\Bx,\BF)$, given by \rf{Wcirc}. Observe,
that the Morrey quasiconvexity condition \rf{qxi} can be written as
\beq
\int_{B(\Bzr,1)} W^{\circ}(\Bx_{0}, \Grad\BGf (\Bx) )d\Bx\geq 0
\eeq{qxipr}
for all $\BGf\in W_{0}^{1,\infty}(B(\Bzr,1);\bb{R}^{m})$, because, clearly
\[
\int_{B(\Bzr,1)}(W_{\BF}(\Bx_{0},\BF(\Bx_{0})),\Grad\BGf(\Bx))d\Bx=0.
\]
If $d=1$ or $m=1$, condition \rf{qxipr} reduces to the Weierstrass condition
$W^{\circ}(\Bx,\BF)\ge 0$ for all $\Bx$ and $\BF$.
Similarly to \rf{qxipr}, quasiconvexity at the free boundary condition
\rf{qxb} can be written as 
\beq
\int_{B^{-}_{\Bn(\Bx_{0})}(\Bzr,1)}W^{\circ}(\Bx_{0},\Grad\BGf (\Bx) )d\Bx\geq 0
\eeq{qxbpr}
for all $\BGf\in W_{0}^{1,\infty}(B(\Bzr,1);\bb{R}^{m})$, because
\beq
\int_{B^{-}_{\Bn(\Bx_{0})}(\Bzr,1)}(W_{\BF}(\Bx_{0},\BF(\Bx_{0})),\Grad\BGf(\Bx))d\Bx=0.
\eeq{free}
The vanishing of the integral in \rf{free} occurs because of the boundary
condition in \rf{e-l eqn}. We, however, will regard inequalities \rf{qxipr}
and \rf{qxbpr} as 
\emph{primary} conditions that reduce to \rf{qxi} and \rf{qxb} in case $\By(\Bx)$
satisfies the Euler-Lagrange equation. (Of course, \rf{qxi} and \rf{qxipr} are
unconditionally equivalent.) 

We summarize our discussion of necessary conditions for $W^{1,\infty}$ weak-* 
local minima above in the form of a theorem for reference
purposes.
\begin{theorem}[Necessary conditions]
\label{necessary}
Let $\By\in\CA$ be a $W^{1,\infty}$ weak-*  local minimizer then
\begin{enumerate}
\item[(i)] $\By(\Bx)$ is a weak solution of the Euler-Lagrange equation \rf{e-l eqn}.
\item[(ii)] The second variation \rf{secvar}
is nonnegative for all $\BGf\in\Var(\CA)$.
\item[(iii)] Quasiconvexity inequalities \rf{qxipr} and \rf{qxbpr} hold for
  all $\BGf\in W_{0}^{1,\infty}(B(\Bzr,1);\bb{R}^{m})$.
\end{enumerate}
\end{theorem}

The following theorem corresponds to Step 2 in our ``Strategy'' on page
\pageref{strat} and is the basis for the sufficient conditions for
$W^{1,\infty}$ weak-*  local minima. 
\begin{theorem}
\label{semisuff1}
Let $\By \in C^{1}(\bra{\GO};\bb{R}^{m})$  satisfy conditions (ii) and (iii) 
of Theorem~\ref{necessary}. Then $\Gd'E\ge 0$ for any variation 
$\{\BGf_{n}:n\ge 1\}\subset\Var(\CA)$.
\end{theorem}
\begin{cor} 
\label{newsemisuff1}
Let $\By \in C^{1}(\bra{\GO};\bb{R}^{m})\cap\CA$  satisfy conditions (i)--(iii) 
of Theorem~\ref{necessary}. Then $\Gd E\ge 0$ for any variation 
$\{\BGf_{n}:n\ge 1\}\subset\Var(\CA)$.
\end{cor}
The theorem says that on the size scale determined by
$\|\Grad\BGf_{n}\|_{2}^{2}$ the variation $\{\BGf_{n}\}$ cannot decrease the
value of the
functional. In order to resolve the question of $W^{1,\infty}$ weak-* local minima,
one needs to understand the set of variations resulting in $\Gd' E=0$. We will
call such variations ``neutral''. At the
moment it is still an open problem to characterize all neutral
variations, but, as we show in Theorem~\ref{semisuff2},
a natural strengthening of necessary conditions (ii)--(iii) in
Theorem~\ref{necessary} will be sufficient to eliminate all neutral
variations altogether. We remark, however, that in the presence of jump
discontinuities of $\BF(\Bx)$ the set of neutral variations is never
empty \cite{grtr08}.
Hence, without our assumption of continuity of the gradient $\BF(\Bx)$ the
sufficient conditions in Theorem~\ref{semisuff2} below cannot possibly
be satisfied.
\begin{theorem}[Sufficient conditions]
\label{semisuff2}
Let  $\By \in C^{1}(\bra{\GO};\bb{R}^{m})\cap\CA$ solve the Euler-Lagrange equation
\rf{e-l eqn} weakly. Assume that there exists $\Gb > 0$ such that 
\begin{enumerate}
\item[(ii)'] The second variation is uniformly positive
\[
\Gd^{2}E=\int_{\GO}(\SFL(\Bx)\Grad\BGf(\Bx),\Grad\BGf(\Bx))d\Bx\ge
\Gb\int_{\GO}|\Grad\BGf (\Bx)|^{2}d\Bx
\]
for all $\BGf\in\Var(\CA)$.
\item[(iii)'](Uniform quasiconvexity)
\begin{enumerate}
\item[(a)] for all $\Bx _{0} \in \GO$
\beq
\int_{B(\Bzr,1)} W^{\circ}(\Bx_{0}, \Grad\BGf (\Bx) )d\Bx\geq
\Gb\int_{B(\Bzr,1)}|\Grad\BGf (\Bx)|^{2}d\Bx
\eeq{uqcx}
for all $\BGf\in W_{0}^{1,\infty}(B(\Bzr,1);\bb{R}^{m})$.
\item[(b)] for all $\Bx _{0}\in \partial \GO _{2}$ 
\beq
\int_{B^{-}_{\Bn(\Bx_{0})}(\Bzr,1)}W^{\circ}(\Bx_{0},\Grad\BGf (\Bx) )d\Bx\geq
\Gb\int_{B^{-}_{\Bn(\Bx_{0})}(\Bzr,1)}|\Grad\BGf (\Bx)|^{2}d\Bx
\eeq{uqcxb}
for all $\BGf\in W_{0}^{1,\infty}(B(\Bzr,1);\bb{R}^{m})$.
\end{enumerate}
\end{enumerate}
Then $\Gd E \geq \Gb$ for any variation $\{\BGf_{n}\}$. 
In particular $\By(\Bx)$ is a $W^{1,\infty}$ weak-* local minimizer of $E$.
\end{theorem}
Theorem~\ref{semisuff2} is an immediate corollary of
Theorem~\ref{semisuff1}, as shown in the following proof.
\begin{proof}
Let
\[
\Tld{W}(\Bx,\BF) = W(\Bx,\BF) -\Gb|\BF|^{2}.
\] 
Then 
\[
\Tld{W}^{\circ}(\Bx,\BF) = W^{\circ}(\Bx,\BF) -\Gb|\BF|^{2}.
\] 
Observe that conditions (ii)', (iii)'(a) and (iii)'(b) can be rewritten as
conditions (ii) and (iii) of Theorem~\ref{necessary} for
$\Tld{W}^{\circ}(\Bx,\BF)$. Then, by Theorem~\ref{semisuff1} applied to
$\Tld{W}$ and $\By(\Bx)$
\[
\Gd'\Tld{E}=\limi_{n\to\infty}\nth{\|\Grad\BGf_{n}\|_{2}^{2}}\int_{\GO}
\Tld{W}^{\circ}(\Bx,\Grad\BGf_{n})d\Bx\ge 0.
\]
But $\Gd'\Tld{E}=\Gd' E-\Gb$. Thus, since $\By(\Bx)$ solves \rf{e-l eqn},
\[
\Gd E=\Gd' E=\Gb+\Gd'\Tld{E}\ge\Gb>0.
\]
It follows that for every variation $\{\BGf_{n}\}$ the functional
increment $\GD E(\BGf_{n})$ is non-negative for all $n$ large enough,
and so $\By(\Bx)$ is a $W^{1,\infty}$ weak-* local minimizer.
\end{proof}

The remaining part of the paper is devoted to the proof of
Theorem~\ref{semisuff1}.  The proof is split into several parts.
All but the last of the parts can be regarded as analytical tools, since they are
independent of the assumptions of Theorem~\ref{semisuff1}. 

In Section~\ref{section:rep} we prove a representation formula that emerges
from our idea to examine the effect of a given variation on a whole space of
Lagrangians $W$. In Sections~\ref{sub:dl} and \ref{sub:orth} we discuss two
related recent developments in Analysis, that concern
the ``oscillations'' and ``concentrations'' behavior of a sequence of
gradients of vector fields. A gradient has a very rigid geometric structure.
The fundamental question is the following: if we permit a sequence of
gradients to be unbounded (in $L^{\infty}$) on a ``small'' set, would we be
able to relax some of that geometric rigidity on the complement of that
``small'' set?  It turns out that the answer is negative. Geometric rigidity
appears to be very robust.  This is established by means of the Decomposition
Lemma \cite{fomupe98,kris94} (see Lemma~\ref{lem:decomp} in
Section~\ref{sub:dl}) and the Orthogonality principle of
Section~\ref{sub:orth} (which we gleaned from one of the technical steps in
\cite{fomupe98}). These two results say that a sequence of gradients that are
unbounded in $L^{\infty}$ (but bounded in $L^{p}$) can be decomposed into
\emph{non-interacting}, or ``orthogonal'' parts, one of which is responsible
only for the oscillations, while the other is responsible only for the
concentrations. At the same time both components retain rigid gradient
structure of the original sequence. The concentration part ``lives'' in some
sense on a set of zero Lebesgue measure,\footnote{We will show by an example
  that this is actually false. However, this image does help on an intuitive
  level.} and can be represented as a ``superposition'' of variations of the
type \rf{nuc}. In order to make the last idea rigorous we adapt the
Localization Principle---a standard technique in the study of Young measures
\cite[Section 8.2]{pedrbk97}---to our setting. The tools developed so far deal
with actions of variations on Lagrangians. As such, they do not require any of
the necessary conditions for local minima to be satisfied. In
Section~\ref{prf-theorem} we combine the tools from the preceeding sections
and the necessary conditions (ii) and (iii) of Theorem~\ref{necessary} to
complete the proof of Theorem~\ref{semisuff1}. We must mention that the same
sequence of steps as in this paper: the representation formula, the decomposition
lemma, the orthogonality principle and the localization principle, was used in
\cite{fomupe98} to characterize the weak-* limits of a
non-linear transformation of the sequence of gradients.


\section{The representation formula}
\setcounter{equation}{0}
\label{section:rep}
\begin{theorem}
\label{existence-main}
Let $\BGy_{n}$  be a bounded sequence in the Sobolev space 
$W^{1,2}(\GO;\bb{R}^{m})$. Suppose 
$\Ga_{n}$ is a sequence of positive numbers such that
$\BGf_{n}(\Bx)=\Ga_{n}\BGy_{n}(\Bx)$ is bounded in $W^{1,\infty}(\GO;\bb{R}^{m})$.
Let 
\beq
R=\sup_{n\ge 1}\|\Grad\BGf_{n}(\Bx)\|_{\infty}.
\eeq{Rbound1}
Then there exist a subsequence, not relabeled, a
nonnegative Radon measure $\Gp$
on $\bra{\GO}$, and  families of  probability measures 
$\{\Gm_{\Bx}\}_{\Bx\in\bra{\GO}}$ supported on the ball $\bra{\CB(\Bzr,R)}$ in
$\bb{R}^{m\times d}$ and $\{\Gl_{\Bx}\}_{\Bx\in\bra{\GO}}$ supported on the
unit sphere $\CS$ in $\bb{R}^{m\times d}$ with the property that
\beq
\CF(\Bx, \Ga_{n}, \Grad\BGy_{n}) \wk*\CI (\Bx, \mu_{\Bx}, \Gl_{\Bx})d\Gp
\eeq{mdef1}
in the sense of measures, where $\CF(\Bx,\Ga,\BG)$ is given by \rf{shhd} and
\beq
\CI(\Bx, \mu_{\Bx}, \Gl_{\Bx}) = \int_{\bra{\CB(\Bzr,R)}}U(\Bx,
\BF)d\mu_{\Bx}(\BF) + \hf\int_{\CS}(\SFL (\Bx)\BF, \BF)d\Gl_{\Bx}(\BF).
\eeq{Idef}
In particular, $|\Grad \BGy_{n}|^{2} \wk* d\Gp$.
  \end{theorem}
Note that in this theorem we do not assume that $\Ga_{n}=\|\Grad\BGf_{n}\|_{2}$. 
\begin{proof}
For each $n\ge 1$, consider a measure $d\Gp_{n} = |\Grad\BGy_{n}(\Bx)|^{2}d\Bx$
on $\bra{\GO}$ and a map $\BGF_{n}:\bra{\GO}\to\bra{\GO}\times \bra{\CB(\Bzr,R)}$,
given by 
\[
\BGF_{n}(\Bx) = (\Bx, \Grad\BGf_{n}(\Bx)).
\]
Let the measure $M_{n}$ on $\bra{\GO}\times \bra{\CB(\Bzr,R)}$ be the 
push-forward of $d\Gp_{n}$ by $\BGF_{n}$. Then, 
for any continuous function $U(\Bx, \BF)$, we have 
\[
\int_{\bra{\GO}\times \bra{\CB(\Bzr,R)}}U(\Bx, \BF) dM _{n}(\Bx, \BF) =
\int_{\bra{\GO}}U(\BGF_{n}(\Bx))d\pi_{n}(\Bx)=
\int_{\GO}U(\Bx, \Grad\BGf_{n}(\Bx))|\Grad \BGy_{n}(\Bx)|^{2}d\Bx.
\] 
From this formula it is clear that $M_{n}$ is a sequence of non-negative
measures and that there exist 
some constant $C > 0$ such that for all $n$, 
$$M_{n} (\bra{\GO}\times\bra{\CB(\Bzr,R)}) = 
\int_{\GO}|\Grad \BGy_{n} (\Bx)|^{2}d\Bx \leq C,$$
since $\Grad \BGy_{n} (\Bx)$ is bounded in $L^{2}$.
That is,  $M_{n}$ is a bounded sequence of measures  in $
\CM(\bra{\GO}\times \bra{\CB(\Bzr,R)})$, where $\CM(\bra{\GO}\times
\bra{\CB(\Bzr,R)})$ is the dual of $C(\bra{\GO}\times \bra{\CB(\Bzr,R)})$.
Then, by the Banach-Alaoglu theorem we can find a
subsequence, not relabeled, and a nonnegative measure  $M$ on 
$\bra{\GO}\times \bra{\CB(\Bzr,R)}$ 
such that $M_{n} \wk* M$ in the sense of measures.
Let $\Gp$ be the projection of $M$ onto $\bra{\GO}$.  Then by the Slicing
Decomposition Lemma \cite{fons92} 
there exists a family of probability measures $\mu=\{\mu_{\Bx}\}_{\Bx\in\bra{\GO}}$ 
on $\bra{\CB(\Bzr,R)}$ such that $M=\mu_{\Bx}\otimes\Gp$ in the sense that
for all $U(\Bx, \BF) \in C(\bra{\GO} \times \bb{R}^{m\times d})$ we have 
\[
\int_{\bra{\GO} \times\bra{\CB(\Bzr,R)}}U(\Bx, \BF)dM (\Bx, \BF) =
\int_{\bra{\GO}}\int_{\bra{\CB(\Bzr,R)}}U(\Bx, \BF)d\mu_{\Bx}(\BF)d\Gp(\Bx)
\]
Therefore, 
\beq
\lim_{n\to\infty}\int_{\GO}
U(\Bx, \Ga_{n}\Grad \BGy_{n}(\Bx))|\Grad\BGy_{n}(\Bx)|^{2}d\Bx = 
\int_{\bra{\GO}}\int_{\bra{\CB(\Bzr,R)}}U(\Bx, \BF)d\mu_{\Bx}(\BF) d\Gp(\Bx).
\eeq{rep1}
Setting $U(\Bx, \BF) = \Gx(\Bx)$, for $\Gx(\Bx) \in C(\bra{\GO})$
we have 
\[
\lim_{n\to\infty}\int_{\GO}\Gx(\Bx)|\Grad \BGy_{n}(\Bx)|^{2}d\Bx  = \int_{\bra{\GO}}\Gx(\Bx)d\Gp(\Bx),
\]
implying that $\Gp_{n} \wk* \Gp$ in the sense of measures.

Consider now the sequence of vector-valued measures
$d\BGb_{n} = |\Grad\BGy_{n}(\Bx)|\Grad\BGy_{n}(\Bx)d\Bx$ with the
polar decomposition (see \cite{fons92})
$d\BGb_{n}=\Hat{\BGb}_{n}(\Bx)d\Gp_{n}(\Bx) $, where
\[
\Hat{\BGb}_{n}(\Bx)=\frac{\Grad\BGy_{n}(\Bx)}{|\Grad\BGy_{n}(\Bx)|}.
\]
Applying the Varifold limit theorem \cite{fons92}
to $d\BGb_{n}$, we obtain a family of
probability measures $\Gl_{\Bx}$ on the unit sphere $\CS$ in $\bb{R}^{m\times d}$  
such that for any $f \in C(\bra\GO\times\CS)$ 
\beq
f(\Bx, \Hat{\BGb}_{n} (\Bx))d\Gp _{n} \wk* \left[\int_{\CS}f(\Bx, \BF)d\Gl_{\Bx}(\BF)
\right]d\Gp
\eeq{fonseca2}
in the sense of measures.
If we choose $f(\Bx, \BF)=(\SFL(\Bx)\BF,\BF)$, where $\SFL(\Bx)$ is given by \rf{Ldef},
we will obtain, according to \rf{fonseca2}, that 
\[
(\SFL(\Bx)\Grad\BGy_{n}(\Bx),\Grad\BGy_{n}(\Bx)) \wk*
\left[\int_{\CS}(\SFL(\Bx)\BF,\BF)d\Gl_{\Bx}(\BF)\right] d\Gp.
\]
Combining that with \rf{rep1} and recalling \rf{shhd} we obtain \rf{mdef1}.
\end{proof}


\section{The decomposition lemma}
\label{sub:dl}
\setcounter{equation}{0}
The decomposition lemma can be found in \cite{fomupe98,kris94} in great
generality. Here we are going to formulate a slightly more restricted version
but with an extra statement that we need and that is easy to obtain from the proof, but
not from the statement of the Lemma in \cite{fomupe98}. For that reason we will
have to revisit the relevant parts of the proof of the Lemma given in \cite{fomupe98}.
\begin{lemma}[Decomposition Lemma] 
\label{lem:decomp}
Suppose the sequence $\{\BGy _{n}:n\ge 1\}\subset \Var(\CA)$ is bounded in 
$W^{1,2}(\GO;\bb{R}^{m})$.  Then there
exist a subsequence $n(j)$ and sequences $\Bv_{j}$, with mean zero, and $\Bz_{j}$ in
$W^{1,\infty}(\GO;\bb{R}^{m})$ such that $\BGy_{n(j)}=\Bz_{j}+\Bv_{j}$, 
$|\Grad\Bz_{j}|^{2}$ is equiintegrable, $\Bv_{j} \weak\Bzr$ weakly in 
$W^{1,2}(\GO;\bb{R}^{m})$. Moreover there exists a sequence of subsets
$R_{j}$ of $\GO$, such that $|R_{j}|\to 0$ as $j\to\infty$ and
\[
\Bz_{j}(\Bx)=\BGy_{n(j)}(\Bx)\text{ and }\Grad\Bz_{j}(\Bx)=\Grad\BGy_{n(j)}(\Bx)
\text{ for all }\Bx\in\GO\setminus R_{j}.
\] 
In addition, if for some sequence $\Ga_{n}$ of positive numbers
the sequence of functions $\Ga_{n}\Grad \BGy_{n}$ is bounded in 
$L^{\infty}(\GO; \bb{R}^{m\times d})$,
then so are the sequences $\Ga_{n(j)}\Grad\Bz_{j}$ and $\Ga_{n(j)}\Grad\Bv_{j}$.
\end{lemma}
After the proof of the Lemma we will restrict our attention to the subsequence
$n(j)$.  For this reason, the symbols $\Ga_{n}$, $\BGy_{n}$,
$\Bz_{n}$ and $\Bv_{n}$ will refer to $\Ga_{n(j)}$, $\BGy_{n(j)}$, $\Bz_{j}$ and
$\Bv_{j}$ respectively.
\begin{proof}
We split the proof into two parts.
In the first part of the proof we are going to recap the construction of
sequences $\Bz_{n}$ and $\Bv_{n}$ in
\cite{fomupe98}. In the second part we are going to use the details of that
construction to prove the last statement in the Lemma.

\textbf{Part I.}
Recall that we have assumed that $\GO$ is a smooth domain. According to 
\cite[Theorem 7.25]{gitrbk01} there exists an extension operator
\[
X: W^{1,p}(\GO;\bb{R}^{m})\to W^{1,p}(\bb{R}^{d};\bb{R}^{m}),\quad 1\le p\le\infty
\]
and a constant $C>0$ \emph{independent of $p$}, such that for all 
$\BGy\in W^{1,p}(\GO;\bb{R}^{m})$ 
\beq
\|X\BGy\|_{W^{1,p}(\bb{R}^{d};\bb{R}^{m})}\le C\|\BGy\|_{W^{1,p}(\GO;\bb{R}^{m})}.
\eeq{ext}

Let  $\BGy _{n}\in \Var(\CA)$ be a bounded sequence in
$W^{1,2}(\GO;\bb{R}^{m})$. We identify $\BGy _{n}$ with its extension
$X\BGy_{n}$. 
Then the sequence of maximal functions $\{M(\Grad\BGy _{n} )\}$ is bounded in
  $L^{2}(\bb{R}^{d})$ (see \cite[Theorem 1(c), p. 5]{stein70}) and 
the sequence $\Ga_{n}M(\Grad\BGy _{n} )$ is bounded in $L^{\infty}$. 
Let $\Gu = \{\Gu_{\Bx}\}_{\Bx\in\GO}$ be the Young measure generated by a 
subsequence $\{M(\Grad\BGy_{n(k)} )\}$. 
Consider the truncation maps $T_{j}:\bb{R} \to \bb{R}$ given by 
\beqn
T_{j}(s) = \left\{\begin{array}{cl}  
s,& |s|\leq j\\[1ex]
\dfrac{js}{|s|},& |s|>j.
\end{array}
\right.
\eeqn 
For each $j$ the function $T_{j}(s)$ is bounded and therefore, the sequence
$\{|T_{j}(M(\Grad\BGy_{n(k)}))|^{2}:k\ge 1\}$ is equiintegrable.
It follows from \cite[Theorem 6.2]{pedrbk97} 
that for each $j$
\[
|T_{j}(M(\Grad\BGy_{n(k)}))|^{2}\weak\int_{\bb{R}}|T_{j}(s)|^{2}d\Gu_{\Bx}(s),
\text{ as }k\to\infty
\]
weakly in $L^{1}(\GO)$.
Let \[ \skew{6}\bra{f}(\Bx) = \int_{\bb{R}}|s|^{2}d\Gu_{\Bx}(s). \]  Then, according
to the theory of Young measures \cite[Theorem 6.11]{pedrbk97},
$\skew{6}\bra{f}\in L^{1}(\GO)$. Notice that $|T_{j}(s)|\le |s|$. Therefore, by the
dominated convergence theorem, we have 
\[
\int_{\bb{R}}|T_{j}(s)|^{2}d\Gu_{\Bx}(s)\weak\skew{6}\bra{f}(\Bx)\text{ as }j\to\infty
\]
weakly in $L^{1}(\GO)$. It turns out that it is possible to choose a subsequence $k(j)$
such that 
\beq
|T_{j}(M(\Grad\BGy_{n(k(j))}))|^{2}\weak\skew{6}\bra{f}(\Bx)\text{ as }j\to\infty
\eeq{dblim}
weakly in $L^{1}(\GO)$ (the proof is given in \cite{fomupe98}\footnote{Since
the space $L^{\infty}$ is not separable, we cannot claim a priori that a
limit of limit points of the sequence is a limit point of the sequence in a
weak topology of $L^{1}$.}).
To simplify notation, let $n(j)$ denote $n(k(j))$.
Set 
\[R'_{j} = \{\Bx\in \GO: M(\Grad\BGy _{n(j)})(\Bx) \geq j\}.\]
Since $\GO$ is bounded and $M(\Grad\BGy _{n(j)})$ is bounded in $L^{2}(\GO)$,
we have $|R'_{j}| \to 0$ as $j\to \infty$.
In \cite[p. 255, Claim \#2]{evga92} it is proved that there exist Lipschitz functions 
$\Bz_{j}'$ such that 
\[\Bz_{j}'(\Bx) = \BGy _{n(j)}(\Bx) \quad \textrm{for a.e. }\Bx\in\bb{R}^{d}\setminus
R'_{j}, \text{ and }|\Grad\Bz_{j}'(\Bx)| \leq Cj~~ \textrm{ for a.e. }\Bx \in
\bb{R}^{d}.\] 
Let 
\[
R_{j} = R'_{j}\cup\{\Bx\in \GO: \Grad\Bz '_{j}(\Bx)\not= \Grad\BGy_{n(j)}(\Bx)\}.
\] 
The sets $R_{j}$ and $R'_{j}$ differ by a set of
Lebesgue measure zero by \cite[Theorem~3 and Remark~(ii), Section~6.1.3]{evga92}. 
Therefore, 
$|R_{j}| \to 0$ as $j\to \infty$.

\textbf{Part II.}
Observe that on $\GO\setminus R_{j}$ we have the inequality
\[
|\Grad z'_{j}(\Bx)|=|\Grad \BGy _{n(j)}(\Bx)|\leq 
|M(\Grad\BGy _{n(j)})(\Bx)| = |T_{j}(M(\Grad\BGy _{n(j)})(\Bx))|
\]
while if $\Bx \in R_{j}'$, then 
\[
|\Grad z'_{j}(\Bx) |\leq C j = C |T_{j}(M(\Grad\BGy _{n(j)})(\Bx))|
\]
We conclude that 
\beq
|\Grad z'_{j}(\Bx) |\leq C |T_{j}(M(\Grad\BGy _{n(j)})(\Bx))|~~\textrm{for
  a.e. }\Bx \in \GO,
\eeq{zcontr}
which, together with \rf{dblim}, yields the equiintegrability of  $\{|\Grad
\Bz'_{j}|^{2}\}$ and boundedness of $\Ga_{n(j)}\Grad\Bz'_{j}$ in  $L^{\infty}$. 

Let, $\Bv'_{j} = \BGy_{n(j)} - \Bz_{j}'$. Then, $\Grad\Bv'_{j}$ is bounded in
$L^{2}$ because so are $\Grad\BGy_{n(j)}$ and $\Grad\Bz'_{j}$ 
(as $|\Grad\Bz'_{j}|^{2}$ is equiintegrable). Similarly,
$\Ga_{n(j)}\Grad\Bv'_{j}$ is bounded in $L^{\infty}$, because so are
$\Ga_{n(j)}\Grad\BGy_{n(j)}$ and $\Ga_{n(j)}\Grad\Bz'_{j}$. 
Now, let $\av{\Bv'_{j}}$ be the average of the field $\Bv'_{j}$ over $\GO$ and let
\[
\Bz_{j} = \Bz '_{j} + \av{\Bv '_{j}},\qquad \Bv_{j} = \Bv '_{j} -
\av{\Bv '_{j}}.
\]
Then, by Poincar\'e inequality, $\Bv_{j}$ is bounded in
$W^{1,2}(\GO;\bb{R}^{m})$. Thus, $\Bz_{j}$
is also bounded in $W^{1,2}$, since $\BGy_{n}$ is. 
Finally, for any $\BGvf\in W^{1,2}(\GO;\bb{R}^{m})$ we have
\beqn
\left|\int_{\GO}(\BGvf, \Bv_{j}(\Bx))d\Bx\right.  &+& \left.\int_{\GO}(\Grad \BGvf
  (\Bx),\Grad\Bv_{j}(\Bx))d\Bx\right|\\
&&\leq \left(\int_{R_{j}}|\BGvf
(\Bx)|^{2}d\Bx\right )^{1/2}\|\Bv_{j}\|_{L^{2}}  +
\left(\int_{R_{j}}|\Grad\BGvf(\Bx)|^{2}
  d\Bx\right)^{1/2}\|\Grad\Bv_{j}\|_{L^{2}} \to 0 
\eeqn
as $j\to \infty$ 
since the sequence $\Bv_{j}$ is bounded in $W^{1,2}(\GO;\bb{R}^{m})$, and
$|R_{j}|\to 0$. This proves that
$\Bv_{j} \weak 0$ in $W^{1,2}(\GO;\bb{R}^{m})$.   
\end{proof}


\section{The orthogonality principle}
\label{sub:orth}
\setcounter{equation}{0}
The decomposition lemma allows us to represent a sequence of gradients that
are bounded in $L^{2}$ as a sum of two sequences of gradients. One of them is
square-equiintegrable and generates the same Young measure as the original sequence,
while the other sequence captures the ``concentration effect''. We are going
to apply the decomposition lemma not to the variation $\{\BGf_{n}\}$ itself
but to the rescaled sequence $\BGy_{n}$ given by \rf{alpsi}. If $\Ga_{n}\to 0$
then the intuitive interpretation of the induced decomposition of $\BGf_{n}$
will be a decomposition of $\BGf_{n}$ into strong ($\Ga_{n}\Bv_{n}$) and weak
variations ($\Ga_{n}\Bz_{n}$), even if the Definition~\ref{def:wkstr} is not
exactly satisfied.

The orthogonality principle says that the two terms in the decomposition of a
variation do not interact (are ``orthogonal''). A version of this lemma was
used in \cite{fomupe98} as one of the steps in their characterization of the
weak-* limits of  of sequences non-linear transformations of gradients.
\begin{lemma}[Orthogonality Principle]
  \label{lem:orth}
\beq
\CF(\Bx,\Ga_{n},\Grad\BGy_{n})-\CF(\Bx,\Ga_{n},\Grad\Bz_{n})-
\CF(\Bx,\Ga_{n},\Grad\Bv_{n})\to 0
\eeq{orth}
strongly in $L^{1}$.
\end{lemma}
Before we prove this lemma, let us combine it with Theorem~\ref{existence-main}.
According to Theorem~\ref{existence-main} in
Section~\ref{section:rep}, there exist measures
$\Tld{M}=\Tld{\mu}_{\Bx}\otimes\Tld{\pi}$ and
$\Tld{\GL}=\Tld{\Gl}_{\Bx}\otimes\Tld{\pi}$ such that
\beq
\CF(\Bx,\Ga_{\Bn},\Grad \Bv_{n}) \wk* 
\CI (\Bx, \Tld{\mu}_{\Bx}, \Tld{\Gl}_{\Bx})d\Tld{\Gp},
\eeq{concentration}
where the functional $\CI$ is given by \rf{Idef}.

We can actually say more about the term involving $\Bz_{n}$ in \rf{orth}.
Let $\nu= \{\nu_{\Bx}\}_{\Bx \in \GO}$ be the gradient Young measure generated
by the sequence $\{\Grad\BGy _{n}\}$. Observe that the sequence
$\{\Grad\Bz_{n}\}$ generates the same Young measure as $\{\Grad\BGy _{n}\}$ because 
$\Grad \Bz _{n} (\Bx) = \Grad \BGy_{n}(\Bx)$ for $\Bx\not\in R_{n}$ and $|R_{n}|
\to 0$ (see \cite[Lemma 6.3(i)]{pedrbk97}).
Moreover since $|\Grad \Bz _{n}|^{2}$ is
equiintegrable, 
\beq
|\Grad \Bz _{n}|^{2} \weak m(\Bx)=\int_{\bb{R}^{m\times d}}|\BF|^{2}d\Gv_{\Bx}(\BF)
\eeq{ym2m}
weakly in $L^{1}(\GO)$.

\begin{lemma}
\label{zeropart}
Assume that $\Ga_{n}\to 0$. Then there exists a subsequence (not relabeled)
such that
\[
\CF(\Bx,\Ga_{n},\Grad\Bz_{n}) \weak U(\Bx, \Bzr)m(\Bx) + \hf\int_{\bb{R}^{m\times d}}
(\SFL(\Bx)\BF,\BF)d\nu_{\Bx}(\BF)
\]
weakly in $L^{1}(\GO)$.
\end{lemma}
By construction, $U(\Bx,\Bzr)=0$. We have included this term in
Lemma~\ref{zeropart} in order to emphasize that ``for practical purposes'' the
values of the sequence $\Ga_{n}\Grad\Bz_{n}(\Bx)$ are uniformly small,
justifying our intuitive understanding of $\Ga_{n}\Bz_{n}(\Bx)$ as the ``weak
part'' of the variation $\BGf_{n}$. Furthermore,
we see that the effect of the variation $\Ga_{n}\Bz_{n}$ on the functional can be
described by a quantity that has an intimate relation to the second variation
\rf{secvar}. This relation will be made absolutely precise in
Section~\ref{prf-theorem} by means of \cite[Lemma 8.3]{pedrbk97}.

Using Lemma~\ref{zeropart}, \rf{mdef1} and
\rf{concentration} we can pass to the limit in \rf{orth} to obtain the decomposition
\begin{equation}\label{decomp}
\CI (\Bx, \mu_{\Bx}, \Gl_{\Bx})d\Gp  = \CI (\Bx,\Tld{\mu}_{\Bx},\Tld{\Gl}_{\Bx} )d\Tld{\Gp} + \CY(\Bx)d\Bx,
\end{equation}
in the sense of measures, where
\[
\CY (\Bx) =  \hf\int_{\bb{R}^{m\times d}}(\SFL(\Bx)\BF,\BF)d\nu_{\Bx}(\BF).
\]

The representation \rf{decomp} holds for any continuous function $U(\Bx,\BF)$
on $\bra{\GO}\times\bb{R}^{m\times d}$ and any continuous fourth order
tensor $\SFL(\Bx)$ on $\bra{\GO}$.  Thus
taking $U = 0$ and $\SFL (\Bx) = \SFI$, the fourth order identity tensor, in
\rf{decomp} we get the decomposition
\beq
d\Gp  = d\Tld{\Gp} + \hf m(\Bx)d\Bx,
\eeq{pidecomp}
where $m(\Bx)$ is defined in \rf{ym2m}. The first term is generated by a
sequence $|\Grad\Bv_{n}|^{2}$ which is non-zero on the sets $R_{n}$ of
vanishing Lebesgue measure, while $m(\Bx)$ is generated by the equiintegrable
part $|\Grad\Bz_{n}|^{2}$ of $|\Grad\BGy_{n}|^{2}$. It would then be
reasonable to assume that the decomposition \rf{pidecomp} is a Lebesgue
decomposition of the measure $\pi$ into the absolutely continuous and singular
parts. Surprisingly, this is false, as is clear from
the following example that is a modification of the 1D example of Ball
and Murat \cite{bamu89}. 
Consider a sequence of functions 
\[
\BGy_{n}(\Bx) = (f_{n}(x_{1}),0,0)
\]
defined on $\GO = [0,1]^{3}$,
where $f_{n}$ is a continuous function on $[0,1]$, such that
\beq
f_{n}'(x) = \left\{ \begin{array}{ll} \frac{n}{\sqrt{2}},&\textrm{when }
x \in \left[\dfrac{k}{n+ 1} - \nth{n^{3}}, \dfrac{k}{n+ 1} +\nth{n^{3}}\right]
\text{ for } k = 0, 1, 2,\ldots, n\\[2ex]
0, & \textrm{otherwise} \end{array}\right.
    \eeq{example}

Then  $|\Grad \BGy_{n}|^{2} = (f_{n}'(x_{1}))^{2} \wk* d\Bx$ in the sense of
measures.  Moreover the Young measure generated by $\Grad \BGy_{n}$ is
$\delta_{\Bzr}$, 
and so $m(\Bx)  = 0$ for all $\Bx \in \GO$.

We conclude this section with proofs of Lemmas~\ref{lem:orth} and \ref{zeropart}.
\begin{proof}[Proof of Lemma~\ref{lem:orth}]
~\\[-2ex]

\textbf{Step 1.} Let's write 
\[ \CF(\Bx,\Ga_{n},\Grad\BGy_{n})-\CF(\Bx,\Ga_{n},\Grad\Bz_{n})-
\CF(\Bx,\Ga_{n},\Grad\Bv_{n})= I_{n}(\Bx; U) + J_{n}(\Bx),
\]
where 
\[ 
I_{n}(\Bx; U)= U(\Bx,\Ga_{n}\Grad\BGy_{n})|\Grad\BGy_{n}|^{2} - 
U(\Bx,\Ga_{n}\Grad\Bv_{n})|\Grad\Bv_{n}|^{2}-
U(\Bx,\Ga_{n}\Grad\Bz_{n})|\Grad\Bz_{n}|^{2}
\] 
and
\[
2J_{n}(\Bx) = (\SFL (\Bx)\Grad\BGy_{n} (\Bx),\Grad\BGy_{n}(\Bx))- 
(\SFL (\Bx)\Grad\Bv_{n} (\Bx),\Grad\Bv_{n}(\Bx))- 
(\SFL (\Bx)\Grad\Bz_{n} (\Bx),\Grad\Bz_{n}(\Bx)).
\]
Therefore to prove the lemma it suffices to show that $I_{n}(\Bx; U) \to 0$
and $J_{n}(\Bx)\to 0$ strongly in $L^{1}$. 

\textbf{Step 2.} Assume that $U$ is smooth. Let us show that $I_{n}(\Bx; U) \to
0$ strongly in $L^{1}$ as $n\to \infty$.
We have
\[
\begin{split}
\int_{\GO}|I_{n}(\Bx; U)|d\Bx \leq  
\int_{R_{n}}\left|U(\Bx,\Ga_{n}\Grad\BGy_{n}(\Bx))|\Grad\BGy_{n}(\Bx)|^{2}\right. &-
  \left.U(\Bx,\Ga_{n}\Grad\Bv_{n}(\Bx))|\Grad\Bv_{n}(\Bx)|^{2}\right|d\Bx\\ 
& + \int_{R_{n}}|U(\Bx,\Ga_{n}\Grad\Bz_{n}(\Bx))||\Grad\Bz_{n}(\Bx)|^{2}d\Bx.
\end{split}
\]
Let 
\beq
R = \sup_{n\ge 1}\left\{\max
(\|\Ga_{n}\Grad \BGy_{n}\|_{\infty},\|\Ga_{n}\Grad \Bz_{n}\|_{\infty},
\|\Ga_{n}\Grad \Bv_{n}\|_{\infty})\right \}.
\eeq{Rbound2}
By mean value theorem, there exists $C = C(R) > 0$
such that 
\beq
|U(\Bx,\BA)|\BA|^{2}- U(\Bx,\BB)|\BB|^{2}| \leq C(|\BA| + |\BB|)|\BA-\BB|
\eeq{uniformestimate} 
for every $\Bx\in \bra\GO$, $|\BA|\le R$ and $|\BB|\le R$. 
Thus we have  
\[
\int_{\GO}|I_{n}(\Bx; U)|d\Bx\leq C\int_{R_{n}}\{
|\Grad\BGy_{n}(\Bx)||\Grad\Bz_{n}(\Bx)|+|\Grad\Bv_{n}(\Bx)||\Grad\Bz_{n}(\Bx)|
+|\Grad\Bz_{n}(\Bx)|^{2}\}d\Bx.
\]
Applying the Cauchy Schwartz inequality to the first two summands on the right hand
side of the  above inequality we get 
\[
\begin{split}
\int_{\GO}|I_{n}(\Bx; U)|d\Bx&\leq
C\|\Grad\BGy_{n}\|_{L^{2}}\left(\int_{R_{n}}|\Grad\Bz_{n}(\Bx)|^{2}d\Bx\right)
^{1/2}\\ &+
C\|\Grad\Bv_{n}\|_{L^{2}}\left(\int_{R_{n}}|\Grad\Bz_{n}(\Bx)|^{2}d\Bx\right)
^{1/2}
+ C\int_{R_{n}}|\Grad\Bz_{n}(\Bx)|^{2}d\Bx.
\end{split}
\]
Equiintegrability of $\Bz_{n}$ and $L^{2}$ boundedness of $\Grad\BGy_{n}$ and
$\Grad\Bv_{n}$ implies that $\|I_{n}(\Bx; U)\|_{1} \to 0$.  

\textbf{Step 3.} Here we show $I_{n}(\Bx; U)\to 0$ strongly in $L^{1}$ as
$n\to \infty$ for all $U$ continuous. Let us approximate  $U$ by
a smooth function. For $\Ge >0$ 
there exists a smooth function $V$ such
that $\|U-V\|_{\infty} < \Ge$ on $\bra\GO \times \bra{\CB(\Bzr,R)}$ .
Then $I_{n}(\Bx;U) =I_{n}(\Bx;V) + I_{n}(\Bx;U - V) $  
and 
\[
\int _{\GO}|I_{n}(\Bx; U-V)|d\Bx \leq
\|U-V\|_{\infty}\left( \|\Grad\BGy_{n}\|_{2}^{2} +
  \|\Grad\Bv_{n}\|_{2}^{2} +\|\Grad\Bz_{n}\|_{2}^{2} \right).
\]
Thus, we get the inequality
\[
\|I_{n}(\Bx;U)\|_{1}\le\|I_{n}(\Bx;V)\|_{1}+C\|U-V\|_{\infty},
\]
from which it follows, by way of Step 1, that $\|I_{n}(\Bx;U)\|_{1}\to 0$.

\textbf{Step 4.} 
The decomposition $\BGy_{n} = \Bz_{n} + \Bv_{n}$ gives 
\beqn
J_{n}(\Bx) = (\SFL (\Bx)\Grad\Bz_{n} (\Bx),\Grad\Bv_{n}(\Bx)). 
\eeqn
It follows that
\[
\int_{\GO}|J_{n}(\Bx)|d\Bx 
\leq C \int_{R_{n}}|\Grad\Bv_{n}(\Bx)||\Grad\Bz_{n}(\Bx)|d\Bx
\leq C\|\Grad\Bv_{n}\|_{2}\left(\int_{R_{n}}|\Grad\Bz_{n}(\Bx)|^{2}\right)^{1/2}\to 0
\]
by the Cauchy Schwartz inequality and the equiintegrability of $|\Grad\Bz_{n}|^{2}$.
This completes the proof of the Lemma.
\end{proof}

\begin{proof}[Proof of Lemma \ref{zeropart}] 
It suffices to prove that 
\beq
U(\Bx,\Ga_{n}\Grad\Bz_{n}(\Bx))|\Grad\Bz_{n}(\Bx)|^{2}\weak U(\Bx,\Bzr)m(\Bx)
\eeq{fass} 
and 
\beq
\hf(\SFL(\Bx)\Grad\Bz_{n}(\Bx),\Grad\Bz_{n}(\Bx)) \weak \CY(\Bx)
\eeq{sass} 
weakly in $L^{1}(\GO)$.
The relation \rf{sass} follows directly from standard theory of Young measures
\cite[Theorem 6.2]{pedrbk97}.
In order to prove \rf{fass} we show that
\beq 
T_{n}(\Bx) = (U(\Bx,\Ga_{n}\Grad\Bz_{n}(\Bx))-
  U(\Bx,\Bzr))|\Grad\Bz_{n}(\Bx)|^{2} \to 0 
\eeq{estim}
strongly in $L^{1}(\GO)$. Then \rf{estim} and the fact that
$|\Grad\Bz_{n}(\Bx)|^{2}\weak m(\Bx)$ weakly in $L^{1}(\GO)$ imply the Lemma.

Let us prove \rf{estim} now.
Observe that $\Ga_{n}\Grad\Bz_{n}\to \Bzr$ in $L^{2}$, because $\Grad \Bz_{n}$
is bounded in $L^{2}$ and $\Ga_{n}\to 0$. Then we can find a subsequence, not
 relabeled,  such that $\Ga_{n}\Grad\Bz_{n}(\Bx) \to \Bzr$ for a.e. $\Bx\in
 \GO$. Let us fix $\Ge > 0$. Then, by the equiintegrability of
 $|\Grad\Bz_{n}|^{2}$, there exists $\Gd >0$ such that
\beq
\sup_{n\ge 1}\int_{E}|\Grad\Bz_{n}|^{2}d\Bx < \Ge,
\eeq{equii}
whenever  $E$ is measurable
and $|E|<\Gd$. Applying Egorov's theorem, we can find the set $E\subset\GO$, 
such that $|E| <\Gd$ and
$\Ga_{n}\Grad\Bz_{n}(\Bx) \to \Bzr$ uniformly on $\GO\setminus E$.
By continuity of $U$, we can 
 find $N\ge 1$ such that for all $n\geq N$ and for all $\Bx\in \GO\setminus E$
we have $|U(\Bx, \Ga_{n}\Grad\Bz_{n}(\Bx)) - U(\Bx, \Bzr)|\leq \Ge$.
At the same time we have  $|U(\Bx, \Ga_{n}\Grad\Bz_{n}(\Bx)) - U(\Bx,
\Bzr)|\leq C$ for all $\Bx \in \GO$, since
$\Ga_{n}\Grad\Bz_{n}(\Bx)$ is bounded in $L^{\infty}$.
Then for all $n\geq N$ we have
\[
\|T_{n}\|_{1}\le\Ge\int_{\GO\setminus E}|\Grad\Bz_{n}(\Bx)|^{2}d\Bx+
C\int_{E}|\Grad \Bz_{n}(\Bx)|^{2}d\Bx.
\]
Using \rf{equii}, we get
\[
\|T_{n}\|_{1}\le\Ge\|\Grad\Bz_{n}\|_{2}^{2}+C\Ge.
\]
We conclude that $T_{n}\to 0$ in $L^{1}$, since $\Grad\Bz_{n}$ is bounded in
$L^{2}$. This finishes the proof of Lemma~\ref{zeropart}.
\end{proof} 


\section{The localization principle}
\setcounter{equation}{0}
\label{sec:loc}
The orthogonality principle reduces the computation of 
$\int_{\GO}\CF(\Bx,\Ga_{n},\Grad\BGy_{n})d\Bx$ to
the computation of the same quantity for $\Bz_{n}$ and $\Bv_{n}$. We saw
in Section~\ref{sub:orth} that the $\Bz_{n}$ part produces the second
variation of the functional in the same way that weak variations \rf{wkvar} do.
We thus, have a direct link between the requirement of positivity of second
variation \rf{secvar} and the non-negativity of the functional increment
corresponding to the variations $\Ga_{n}\Bz_{n}$ (we will make this precise in
Section~\ref{prf-theorem}).  

As we mentioned at the beginning of Section~\ref{sub:orth}, the variation
$\Ga_{n}\Bv_{n}$ should be regarded intuitively as a ``strong part'' of the
variation $\BGf_{n}$.  For that reason, we expect it to be connected in
some way to the quasiconvexity conditions \rf{qxipr}--\rf{qxbpr}. This,
however, is not so clear.  The basic problem is that the variation
$\Ga_{n}\Bv_{n}$ seems to have a global character,\footnote{Even though
  $\Bv_{n}$ ``lives'' on $R_{n}$ with vanishing Lebesgue measure, we know
  nothing about the geometry of the set $R_{n}$. Example \rf{example} shows
  that the character of the variation $\Ga_{n}\Bv_{n}$ can be global indeed.}  
while the quasiconvexity conditions \rf{qxipr}--\rf{qxbpr} are localized at a
single point.  This is exactly where the 
localization principle comes in. It says that the effect of $\Ga_{n}\Bv_{n}$
can be localized at a single point, providing us with the necessary link to
quasiconvexity conditions. Our localization principle is very similar (on a
technical level) to the localization principle for Young measures
\cite[Theorem 8.4]{pedrbk97}, and both can be regarded as versions of the 
Lebesgue differentiation theorem. In our notation the localization principle
can be stated as 
\beq
\CI (\Bx_{0}, \Tld{\mu}_{\Bx_{0}}, \Tld{\Gl}_{\Bx_{0}})=
\lim_{r\to 0}\lim_{n \to \infty}\nth{\Tld{\Gp}(B_{\GO}(\Bx_{0}, r))}
\int_{B_{\GO}(\Bx_{0}, r)}\CF(\Bx_{0},\Ga_{n},\Grad\Bv_{n})d\Bx
\eeq{local1}
for $\Tld{\Gp}$ a.e. $\Bx_{0}\in\bra{\GO}\cap\text{supp}(\Tld{\Gp})$,
where $B_{\GO}(\Bx_{0},r) = B(\Bx_{0}, r)\cap \bra{\GO}$. The problem with \rf{local1}
is that the maps $\Bv_{n}$ do not necessarily have the proper boundary
conditions to be used as test functions $\BGf$ in the quasiconvexity
inequalities \rf{qxipr} and \rf{qxbpr}. In addition, as far as the
quasiconvexity at the boundary \rf{qxbpr} is concerned, the domain
$B_{\GO}(\Bx_{0}, r)$ (or its rescaled version $B^{-}_{r}=(B_{\GO}(\Bx_{0},r)-\Bx_{0})/r$)
is not quite the domain required in \rf{qxbpr}. In this section we
prove a bit more involved versions of \rf{local1} that remedy the
above stated shortcomings.

\begin{theorem}[Localization principle in the interior]
  \label{th:loci}
Let $\Bx_{0}\in \GO \cup \bra{\partial\GO _{1}}$. Let the  cut-off functions $\Gth_{k}^{r}(\Bx)\in
C^{\infty}_{0}(B_{\GO}(\Bx_{0}, r))$ be such that $\Gth_{k}^{r}(\Bx)\to \Gc _{B_{\GO}(\Bx_{0},
  r)}(\Bx)$, while remaining uniformly bounded in $L^{\infty}$. Let
$\Bv_{n}\weak\Bzr$ weakly in $W^{1,2}(\GO;\bb{R}^{m})$. Let $\Ga_{n}$ be a
sequence of positive numbers such that $\Ga_{n}\Bv_{n}$ is bounded in 
$W^{1,\infty}(\GO;\bb{R}^{m})$. Let
$\Tld{M}=\Tld{\mu}_{\Bx}\otimes\Tld{\pi}$ and
$\Tld{\GL}=\Tld{\Gl}_{\Bx}\otimes\Tld{\pi}$ be the measures corresponding to
the pair $(\Ga_{n},\Bv_{n})$ via Theorem~\ref{existence-main}.
Then for $\Tld{\pi}$ a.e. $\Bx_{0}\in\GO\cup\bra{\dOm_{1}}$
\beq
\lim_{r\to 0}\lim_{k\to \infty}\lim_{n\to\infty}
\nth{\Tld{\Gp}(B_{\GO}(\Bx_{0}, r))}
\int_{B_{\GO}(\Bx_{0}, r)}\CF (\Bx_{0}, \Ga_{n},\Grad(\Gth_{k}^{r}(\Bx)\Bv_{n}(\Bx)))d\Bx
  = \CI (\Bx_{0}, \Tld{\mu}_{\Bx_{0}}, \Tld{\Gl}_{\Bx_{0}}) 
\eeq{localinterior} 
\end{theorem}
In order to formulate the localization principle for the free boundary we have
to take care not only of the \bc s, but also of the geometry of the domain,
that is required to have a ``flat'' part of the boundary with the outer
unit normal $\Bn(\Bx_{0})$. We observe that for
smooth domains $\GO$ the set 
\beq
B^{-}_{r}=\frac{B_{\GO}(\Bx_{0},r)-\Bx_{0}}{r}
\eeq{Bplusr}
is ``almost'' the half-ball $B^{-}_{\Bn(\Bx_{0})}(\Bzr,1)$. As $r\to 0$ the set $B^{-}_{r}$
``converges'' to $B^{-}_{\Bn(\Bx_{0})}(\Bzr,1)$. Formally, we say that there exists a family
of diffeomorphisms
$\Bf_{r}:B^{-}_{\Bn(\Bx_{0})}(\Bzr,1)\to B^{-}_{r}$ such that $\Bf_{r}(\Bx)\to\Bx$ in 
$C^{1}(B^{-}_{\Bn(\Bx_{0})}(\Bzr,1))$ and $\Bf_{r}^{-1}(\By)\to\By$ in 
$C^{1}(B^{-}_{r})$ in the sense that
\[
\sup _{\By \in B^{-}_{r}} |\Bf^{-1}_{r}(\By) - \By|\to 0
\text{ and }
\sup_{\Bx \in B^{-}_{r}} |\Grad \Bf^{-1}_{r}(\By) - \BI| \to 0, \text{ as }r\to 0.
\]
Let
\beq
\Bv_{n}^{r}(\Bx) = \frac{\Bv _{n}(\Bx_{0} + r\Bf_{r}(\Bx)) -\BC_{n}^{r}(\Bx_{0})}{r}
\eeq{blupv}
be the blown-up version of $\Bv_{n}$ defined on $B^{-}_{\Bn(\Bx_{0})}(\Bzr,1)$,
where the constants
\[
\BC_{n}^{r}(\Bx_{0}) =
\nth{|B^{-}_{\Bn(\Bx_{0})}(\Bzr,1)|}\int_{B^{-}_{\Bn(\Bx_{0})}(\Bzr,1)}\Bv_{n}(\Bx_{0} + r \Bf_{r}(\Bx)) d\Bx
\]
are chosen such that $\Bv_{n}^{r}(\Bx)$ has zero mean over
$B^{-}_{\Bn(\Bx_{0})}(\Bzr,1)$. 

\begin{theorem}
  \label{th:locb}
Let $\Bx_{0}\in\dOm_{2}\cap\textrm{supp}(\Tld{\pi})$ and let $\Bv_{n}$ and $\Ga_{n}$ be as
in Theorem~\ref{th:loci}. Let $\Bv_{n}^{r}$ be defined by \rf{blupv} and let
the cut-off functions $\Gth_{k}(\Bx)\in C_{0}^{\infty}(B(\Bzr,1))$ be such
that $\Gth _{k}(\Bx)\to \Gc _{B(\Bzr,1)}(\Bx)$, while remaining uniformly 
bounded in $L^{\infty}$. Let $\BGz_{n,k}^{r}(\Bx)=\Gth_{k}(\Bx)\Bv_{n}^{r}(\Bx)$. 
Then
\beq
\lim_{r\to 0}\lim_{k\to\infty}\lim_{n\to
  \infty}\frac{r^{d}}{\Tld{\Gp}(B_{\GO}(\Bx_{0},r))}\int_{B^{-}_{\Bn(\Bx_{0})}(\Bzr,1)}
\CF(\Bx_{0}, \Ga_{n},\Grad \BGz_{n,k}^{r}(\Bx))d\Bx = 
\CI (\Bx_{0}, \Tld{\mu}_{\Bx_{0}},\Tld{\Gl}_{\Bx_{0}}) 
\eeq{locb}
for $\Tld{\Gp}$-a.e $\Bx_{0} \in \dOm_{2}$.  
\end{theorem}

\subsection{Proof of Theorem~\ref{th:loci}}
\label{sub:pri}
\textbf{Step 1.}
We begin by showing that the gradient of the cut-off functions $\Gth_{k}^{r}$ does
not influence the limit in \rf{localinterior}.
\begin{lemma}
    \label{lem:thetak}
For each fixed $k$ and $r$
\[
\lim_{n\to\infty}\int_{B_{\GO}(\Bx_{0},r)}
\CF(\Bx_{0},\Ga_{n},\Grad(\Gth_{k}^{r}(\Bx)\Bv_{n}(\Bx)))d\Bx =
\lim_{n\to\infty}\int_{B_{\GO}(\Bx_{0}, r)}
\CF(\Bx_{0},\Ga_{n},\Gth_{k}^{r}(\Bx)\Grad\Bv_{n}(\Bx))d\Bx.
\]
  \end{lemma}
\begin{proof}
Let
\[
T_{n,k,r}(\Bx)=\CF(\Bx_{0},\Ga_{n},\Grad(\Gth_{k}^{r}(\Bx)\Bv_{n}(\Bx)))-
\CF(\Bx_{0},\Ga_{n},\Gth_{k}^{r}(\Bx)\Grad\Bv_{n}(\Bx)).
\]
In order to prove the Lemma,
we need to estimate $T_{n,k,r}(\Bx)$ and prove that 
\beq
\int_{B_{\GO}(\Bx_{0},r)}T_{n,k,r}(\Bx)d\Bx\to 0,\text{ as }n\to\infty.
\eeq{errorT}
Notice that our smoothness assumptions on $W$ implies that 
\beq
|\CF (\Bx,\Ga,\BG_{1})-\CF (\Bx,\Ga,\BG_{2})| \leq 
C(M)|\BG_{1} - \BG_{2}|(|\BG_{1}|+|\BG_{2}|)
\eeq{F-estimate}
for some positive constant $C(M)$, when $|\BG_{1}|\leq M$ and $|\BG_{2}|\leq M$.
Therefore,
\[
|T_{n,k,r}(\Bx)|\le C(k,r)|\Grad\Gth_{k}^{r}(\Bx)||\Bv_{n}(\Bx)|(|\Gth_{k}^{r}(\Bx)|
|\Grad\Bv_{n}(\Bx)|+|\Grad\Gth_{k}^{r}(\Bx)||\Bv_{n}(\Bx)|),
\]
which implies that \rf{errorT} holds, because $\Bv_{n}\weak\Bzr$ in $W^{1,2}$.
\end{proof}

\textbf{Step 2.} Next we compute the limit in Lemma~\ref{lem:thetak} by means of
Theorem~\ref{existence-main} and show that the limit in $k\to\infty$
corresponds to taking $\Gth_{k}^{r}(\Bx)=\Gc _{B_{\GO}(\Bx_{0},r)}(\Bx)$.
 \begin{lemma}\label{claim1} 
\beq
\lim_{k\to \infty}\lim_{n\to\infty}
\int_{B_{\GO}(\Bx_{0}, r)}\CF (\Bx_{0}, \Ga_{n},\Gth_{k}^{r}(\Bx)\Grad\Bv_{n}(\Bx))d\Bx 
 = \int_{B_{\GO}(\Bx_{0},r)} \Tld{\CI}(\Bx_{0},\Bx)d\Tld{\Gp}(\Bx)
\eeq{intermediate}
where 
\beq
\Tld{\CI}(\Bx_{0}, \Bx)  =  \int_{\bra{\CB(\Bzr,R)}} U(\Bx_{0},
  \BF)d\Tld{\mu}_{\Bx}(\BF) + \hf\int_{\CS} (\SFL (\Bx_{0})\BF, \BF)d\Tld{\Gl}_{\Bx}(\BF),
\eeq{Itld}
where $R$ is given by \rf{Rbound2}.
\end{lemma}
\begin{proof}
For each fixed $\Bx_{0}\in\GO\cup\bra{\dOm_{1}}$ and $k\ge 1$ we define
\[
\Tld{U}_{\Bx_{0}}^{(k,r)}(\Bx,\BF)=\Gth_{k}^{r}(\Bx)^{2}U(\Bx_{0},\Gth_{k}^{r}(\Bx)\BF),\quad
\Tld{\SFL}_{\Bx_{0}}^{(k,r)}(\Bx)=\Gth_{k}^{r}(\Bx)^{2}\SFL(\Bx_{0}).
\]
Then, 
\[
\CF(\Bx_{0},
\Ga_{n},\Gth_{k}^{r}(\Bx)\Grad\Bv_{n})=\Tld{\CF}(\Bx,\Ga_{n},\Grad\Bv_{n}),
\]
where $\Tld{\CF}$ is the functional $\CF$, given by \rf{shhd} with $U$ and
$\SFL$ replaced by $\Tld{U}_{\Bx_{0}}^{(k,r)}$ and $\Tld{\SFL}_{\Bx_{0}}^{(k,r)}$
respectively. Applying Theorem~\ref{existence-main}, we obtain
\[
\begin{split}
\lim_{n\to\infty}&\int_{B_{\GO}(\Bx_{0}, r)}\CF(\Bx_{0}, \Ga_{n},\Gth_{k}^{r}(\Bx)\Grad\Bv_{n}(\Bx))d\Bx\\
&= \int_{B_{\GO}(\Bx_{0}, r)}\Gth_{k}^{r}(\Bx)^{2}\left(\int_{\bra{\CB(\Bzr,R)}}U(\Bx_{0},\Gth_{k}^{r}(\Bx)\BF)
    d\Tld{\Gm}_{\Bx}(\BF) +\hf\int_{\CS}(\SFL(\Bx_{0})\BF,
 \BF)d\Tld{\Gl}_{\Bx}(\BF)\right) d\Tld{\Gp}(\Bx).
\end{split}
\]
By bounded convergence theorem, using the fact that
$\Gth_{k}^{r}(\Bx) \to \Gc_{B_{\GO}(\Bx_{0}, r)}(\Bx)$ we have 
\[
\Gth_{k}^{r}(\Bx)^{2} \left(\int_{\bra{\CB(\Bzr,R)}}U(\Bx_{0},\Gth_{k}^{r}(\Bx)\BF)
    d\Tld{\Gm}_{\Bx}(\BF) +\hf\int_{\CS}(\SFL(\Bx_{0})\BF,
 \BF)d\Tld{\Gl}_{\Bx}(\BF)\right) \to \Tld{\CI}(\Bx_{0},\Bx )\Gc_{B_{\GO}(\Bx_{0}, r)}(\Bx)
\]
as $k\to\infty$ for $\Tld{\Gp}$-a.e $\Bx \in \bra{\GO}$.  The conclusion of
the lemma follows from another application of bounded convergence theorem.
\end{proof}

\textbf{Step 3.}
In order to finish the proof of Theorem~\ref{th:loci} we need to
divide both sides of \rf{intermediate} by 
$\Tld{\Gp}(B_{\GO}(\Bx_{0},r))$ and take the limit as $r\to 0$. The result is a
corollary of the ``vector-valued'' version of the Lebesgue differentiation theorem 
\cite[Corollary 2.9.9]{fedebk69}.
Indeed,  $\Tld{\CI}(\Bx_{0}, \Bx)$ is continuous in $\Bx_{0} \in \bra{\GO}$ for
$\Tld{\Gp}$ a.e. $\Bx \in \bra{\GO}$, and 
\[
\int_{\bra{\GO}}\|\Tld{\CI}(\cdot, \Bx)\|_{C(\bra{\GO})}d\Tld{\Gp}(\Bx) < \infty.
\]
Then for any $\Bx_{0}' \in \bra{\GO}$ and for $\Tld{\Gp}$ a.e. $\Bx_{0}\in\bra{\GO}$, we have 
\[\lim_{r \to 0} 
\nth{\Tld{\Gp}(B_{\GO}(\Bx_{0},r))}\int_{B_{\GO}(\Bx_{0},r)}\Tld{\CI}(\Bx_{0}',\Bx)d\Tld{\Gp}(\Bx)= 
\Tld{\CI}(\Bx_{0}',\Bx_{0}).
\]
Setting $\Bx_{0}' = \Bx_{0}$ we obtain \rf{localinterior}. Theorem~\ref{th:loci} is proved.

\subsection{Proof of Theorem~\ref{th:locb}}
\label{sub:prfb}
The proof basically follows the same sequence of steps as the proof of
Theorem~\ref{th:loci} with the only difference that we have to take care not
only of the cut-off functions $\Gth_{k}$ but also of the small 
deformations $\Bf_{r}$. 

\textbf{Step 1.} As in the proof of Theorem~\ref{th:locb}, we first show that
gradients of the cut-off functions $\Gth_{k}$ do not enter the limit
\rf{locb}.
\begin{lemma}
  \label{lem:thetakb}
\beq
\lim_{n\to \infty}\int_{B^{-}_{\Bn(\Bx_{0})}(\Bzr,1)}\CF(\Bx_{0}, \Ga_{n},\Grad \BGz_{n,k}^{r})d\Bx  =
\lim_{n\to \infty}\int_{B^{-}_{\Bn(\Bx_{0})}(\Bzr,1)}
\CF(\Bx_{0}, \Ga_{n},\Gth_{k}(\Bx)\Grad\Bv_{n}^{r}(\Bx))d\Bx.
\eeq{bndryqcx2}
\end{lemma}
The proof is very similar to the proof of Lemma~\ref{lem:thetak} and is therefore
omitted. The more complex dependence of the integrand on $r$ is irrelevant
at this point because $r$ is fixed here.

\textbf{Step 2.} As in the proof of Theorem~\ref{th:loci} we use
Theorem~\ref{existence-main} to compute the limit as $n\to\infty$ and then
pass to the limit as $\Gth _{k}(\Bx)\to\Gc _{B(\Bzr,1)}(\Bx)$. 
Let us change variables 
\beq
\Bx'=\Bx_{0}+r\Bf_{r}(\Bx)
\eeq{chvar} 
in the \rhs\ in \rf{bndryqcx2}. Solving \rf{chvar} for $\Bx$ we get 
\[
\Bx=\Bt_{r}(\Bx')=\Bf_{r}^{-1}((\Bx'-\Bx_{0})/r).
\]
Then
\[
\begin{split}
\lim_{n\to \infty}\int_{B^{-}_{\Bn(\Bx_{0})}(\Bzr,1)}&
\CF(\Bx_{0}, \Ga_{n},\Gth_{k}(\Bx)\Grad\Bv_{n}^{r}(\Bx))d\Bx=\\
&\lim_{n\to \infty}\int_{B_{\GO}(\Bx_{0},r)}\CF(\Bx_{0},\Ga_{n},
\Gth_{k}(\Bt_{r}(\Bx'))\Grad\Bv_{n}(\Bx')\BJ_{r}(\Bx'))
\frac{J_{r}^{-1}(\Bx')}{r^{d}}d\Bx',
\end{split}
\]
where 
\beq
\BJ_{r}(\Bx')=(\Grad\Bf_{r})(\Bf_{r}^{-1}((\Bx'-\Bx_{0})/r)))
\eeq{Jrdef}
and $J_{r}(\Bx')=\det\BJ_{r}(\Bx')$.
Again, as in the proof of Theorem~\ref{th:loci} we represent the
expression under the integral as the functional $\CF$ constructed with
$\Hat{U}_{\Bx_{0}}^{(k,r)}$ and $\Hat{\SFL}_{\Bx_{0}}^{(k,r)}$ replacing $U$ and
$\SFL$, where
\[
\Hat{U}_{\Bx_{0}}^{(k,r)}(\Bx,\BF)=\frac{\Gth_{k}(\Bt_{r}(\Bx)^{2}}{r^{d}J_{r}(\Bx)}
U(\Bx_{0},\Gth_{k}(\Bt_{r}(\Bx)\BF\BJ_{r}(\Bx)))
\frac{|\BF\BJ_{r}(\Bx)|^{2}}{|\BF|^{2}}
\]
and
\[
(\Hat{\SFL}_{\Bx_{0}}^{(k,r)}(\Bx)\BF,\BF)=
\frac{\Gth_{k}(\Bt_{r}(\Bx))^{2}}{r^{d}J_{r}(\Bx)}
(\SFL(\Bx_{0})\BF\BJ_{r}(\Bx),\BF\BJ_{r}(\Bx)).
\]
We remark, that since $U(\Bx,\BF)$ is continuous and $U(\Bx,\Bzr)=0$, then the
same is true for $\Hat{U}_{\Bx_{0}}^{(k,r)}$. Thus, Theorem~\ref{existence-main}
is applicable and the limit as $n\to\infty$ can be computed. 
The passage to the limit as $k\to\infty$ is no different than the same step in
the proof of Theorem~\ref{th:loci}. Thus, we obtain
\[
\lim_{k\to\infty}\lim_{n\to\infty}
\int_{B^{-}_{\Bn(\Bx_{0})}(\Bzr,1)}\CF(\Bx_{0}, \Ga_{n},\Grad
\BGz_{n,k}^{r})d\Bx=\nth{r^{d}}\int_{B_{\GO}(\Bx_{0},r)}
\frac{\Tld{\CI}_{r}(\Bx_{0},\Bx)}{J_{r}(\Bx)}d\Tld{\pi}(\Bx),
\]
where
\[
\Tld{\CI}_{r}(\Bx_{0},\Bx)=\int_{\bra{\CB(\Bzr,R)}}U(\Bx_{0},\BF\BJ_{r}(\Bx))\frac{|\BF\BJ_{r}(\Bx)|^{2}}{|\BF|^{2}}d\Tld{\mu}_{\Bx}(\BF)+
\hf\int_{\CS}(\SFL(\Bx_{0})\BF\BJ_{r}(\Bx),\BF\BJ_{r}(\Bx))d\Tld{\Gl}_{\Bx}(\BF)
\]

\textbf{Step 3.} On this step, we will show that the deformation $\Bf_{r}$
does not influence the limit as $r\to 0$.
\begin{lemma}
  \label{lem:r}
\[
\lim_{r\to 0}\nth{\Tld{\pi}(B_{\GO}(\Bx_{0},r))}\int_{B_{\GO}(\Bx_{0},r)}
\frac{\Tld{\CI}_{r}(\Bx_{0},\Bx)}{J_{r}(\Bx)}d\Tld{\pi}(\Bx)=
\lim_{r\to 0}\nth{\Tld{\pi}(B_{\GO}(\Bx_{0},r))}\int_{B_{\GO}(\Bx_{0},r)}
\Tld{\CI}(\Bx_{0},\Bx)d\Tld{\pi}(\Bx),
\]
where $\Tld{\CI}(\Bx_{0},\Bx)$ is given by \rf{Itld}.
\end{lemma}
\begin{proof}
Observe that $\BJ_{r}(\Bx)\to\BI$, as $r\to 0$ uniformly in $\Bx\in
B_{\GO}(\Bx_{0},r)$ in the sense that
\beq
\lim_{r\to 0}\sup_{\Bx\in B_{\GO}(\Bx_{0},r)}|\BJ_{r}(\Bx)-\BI|=0.
\eeq{unilim}
Indeed, from \rf{Jrdef} it is easy to see that
\[
\sup_{\Bx\in B_{\GO}(\Bx_{0},r)}|\BJ_{r}(\Bx)-\BI|=\sup_{\Bx\in
  B^{-}_{\Bn(\Bx_{0})}(\Bzr,1)}|\Grad\Bf_{r}(\Bx)-\BI|\to 0,
\]
as $r\to 0$. It follows that $J_{r}(\Bx)\to 1$, as $r\to 0$ uniformly in $\Bx\in
B_{\GO}(\Bx_{0},r)$.
We also have that
\beq
\lim_{r\to 0}\sup_{\Bx\in B_{\GO}(\Bx_{0},r)}
\left|\frac{\Tld{\CI}_{r}(\Bx_{0},\Bx)}{J_{r}(\Bx)}-\Tld{\CI}(\Bx_{0},\Bx)\right|=0,
\eeq{aaa}
due to \rf{unilim} and the fact that the measures $\Tld{\mu}_{\Bx}$ and
$\Tld{\Gl}_{\Bx}$ are supported on compact sets.
The Lemma now follows from \rf{aaa} and the estimate
\[
\nth{\Tld{\pi}(B_{\GO}(\Bx_{0},r))}\int_{B_{\GO}(\Bx_{0},r)}
\left|\frac{\Tld{\CI}_{r}(\Bx_{0},\Bx)}{J_{r}(\Bx)}-\Tld{\CI}(\Bx_{0},\Bx)\right|
d\Tld{\pi}(\Bx)\le\sup_{\Bx\in B_{\GO}(\Bx_{0},r)}
\left|\frac{\Tld{\CI}_{r}(\Bx_{0},\Bx)}{J_{r}(\Bx)}-\Tld{\CI}(\Bx_{0},\Bx)\right|.
\]
\end{proof}

\textbf{Step 4.} The limit in Lemma~\ref{lem:r} is already computed in Step~3
in the proof of Theorem~\ref{th:loci}. This finishes the proof of
Theorem~\ref{th:locb}. 

 
\section{Proof of Theorem \ref{semisuff1}}
\setcounter{equation}{0}                
\label{prf-theorem}
Observe that so far we have been developing analytical \emph{tools}, that is
theorems that do not involve any of the necessary conditions for local minima
listed in Theorem~\ref{necessary}. In this section we will combine the tools
with the inequalities from Theorem~\ref{necessary} to prove
Theorem~\ref{semisuff1}. 

\textbf{Step 1.} First we suppose that the sequence of positive numbers $\Ga_{n}$, 
defined in \rf{alpsi} does not converge to zero (i.e. does not have a
subsequence that converges to zero). Then, 
\[
\Gd'E=\nth{\Ga_{0}^{2}}\int_{\GO}\left(\int_{\bb{R}^{m\times d}}
W(\BF(\Bx)+\BF)d\eta_{\Bx}(\BF)-W(\BF(\Bx))\right)d\Bx,
\]
where $\Ga_{0}$ is a non-zero limit of the sequence $\Ga_{n}$ and $\eta_{\Bx}$
is a Young measure generated by a sequence of gradients $\{\Grad\BGf_{n}\}$
that are bounded in $L^{\infty}$. The term
\[
\int_{\bb{R}^{m\times d}}(W_{\BF}(\BF(\Bx)),\BF)d\eta_{\Bx}(\BF)=0
\]
because the sequence $\Grad\BGf_{n}$ converges to zero in $L^{\infty}$ weak-*.
The non-negativity of $\Gd'E$ now follows from the quasiconvexity assumption
\rf{qxi} and \cite[Theorem 8.14]{pedrbk97}.

\textbf{Step 2.} A more interesting (and complicated) case is when $\Ga_{n}\to 0$.
In this case we have
\beq
\Gd' E = \int_{\bra{\GO}}\CI (\Bx,\mu_{\Bx}, \Gl_{\Bx})d\Gp(\Bx).
\eeq{limit-increment}
and a decomposition \rf{decomp} holds. Thus,
\beq
\Gd' E = \int_{\bra{\GO}}\CI (\Bx,\Tld {\mu}_{\Bx},\Tld{\Gl}_{\Bx})d\Tld{\Gp}(\Bx) 
+ \hf\int_{\GO}\int_{\bb{R}^{m\times d}} (\SFL (\Bx)\BF, \BF)d\nu_{\Bx}(\BF)d\Bx.
\eeq{representation2}
To complete the proof of the Theorem we show that
\beq
\int_{\GO}\int_{\bb{R}^{m\times d}}(\SFL (\Bx)\BF,
\BF)d\nu _{\Bx}(\BF)d\Bx\geq 0
\eeq{osc-ym1}
and 
\beq
\CI (\Bx_{0},
\Tld {\mu}_{\Bx_{0}}, \Tld{\Gl}_{\Bx_{0}}) \geq 0 
\quad\textrm{for $\Tld{\Gp}-$ a.e. $\Bx_{0}\in \bra{\GO}$}.
\eeq{concent-varifold1}

\textbf{Step 3.}  We first prove \rf{osc-ym1}.
Observe that 
 since $\|\Grad \BGy_{n}\|_{2} = 1$ and
$\BGy_{n}|_{\partial {\GO}_{1}} = 0$, there exists $\BGy _{0} \in
W^{1,2}(\GO;\bb{R}^{m})$ satisfying $\BGy_{0}|_{\partial {\GO}_{1}} = 0$ and a
subsequence $\{\BGy_{n}\}$, not relabeled, such that $\BGy_{n} \weak \BGy_{0}$
weakly in $W^{1,2}(\GO;\bb{R}^{m} )$.  Since $\Bv_{n} \weak  0$ weakly in
$W^{1,2}(\GO;\bb{R}^{m} )$, we have $\Bz_{n} \weak \BGy _{0}$ weakly in
$W^{1,2}(\GO;\bb{R}^{m} )$.
By 
\cite[Lemma 8.3]{pedrbk97}, we can find a sequence $\Tld{\Bz}_{n}$ such that
$\Tld{\Bz}_{n} - \BGy _{0} \in W^{1,2}_{0}(\GO;\bb{R}^{m} )$ and $\Grad \Bz _{n}$
and $\Grad \Tld{\Bz} _{n}$ generate the same Young measure $\Gv  = \{\Gv
_{\Bx}\}_{\Bx\in \GO}$. It follows that $\Tld{\Bz}_{n}$ satisfies 
$\Tld{\Bz} _{n}|_{\partial {\GO}_{1}} = 0$. Thus, $\Tld{\Bz}_{n}\in\Var(\CA)$ and
\[
\int_{\GO}(\SFL (\Bx)\Grad \Tld{\Bz} _{n}(\Bx),\Grad \Tld{\Bz}_{n}(\Bx))d\Bx\ge 0
\]
for all $n$, according to the condition (ii) of Theorem~\ref{necessary}.
Taking limit as $n\to \infty$ in the above inequality we obtain \rf{osc-ym1}.

\textbf{Step 4.} On this step we prove the inequality \rf{concent-varifold1}.  
For all $\Bx_{0}\in\GO\cup\bra{\dOm_{1}}$ we have that the functions
$\Gth_{k}(\Bx)\Bv_{n}(\Bx)$ vanish on $\Md B_{\GO}(\Bx_{0},r)$ and therefore,
according to the inequality \rf{qxipr} we have
\[
\int_{B_{\GO}(\Bx_{0}, r)}\CF (\Bx_{0},
\Ga_{n},\Grad(\Gth_{k}(\Bx)\Bv_{n}(\Bx)))d\Bx\ge 0
\]
for all $n$, $k$ and $r$. 
Theorem~\ref{th:loci} then tells us that $\CI (\Bx_{0}, \Tld{\mu}_{\Bx_{0}},
\Tld{\Gl}_{\Bx_{0}})\ge 0$ for $\Tld{\pi}$ almost all $\Bx_{0}\in\GO\cup\bra{\dOm_{1}}$.

For all $\Bx_{0}\in\dOm_{2}$, we use functions $\BGz_{n,k}^{r}(\Bx)$ from the
formulation of Theorem~\ref{th:locb}. These functions
are defined on the half-ball $B^{-}_{\Bn(\Bx_{0})}(\Bzr,1)$ and
vanish on the ``round'' part of the boundary of the half-ball. Therefore,
according to the inequality \rf{qxbpr} we have,
\[
\int_{B^{-}_{\Bn(\Bx_{0})}(\Bzr,1)}\CF(\Bx_{0}, \Ga_{n},\Grad \BGz_{n,k}^{r}(\Bx))d\Bx\ge 0
\]
for all $n$, $k$ and $r$. Theorem~\ref{th:locb} then tells us that 
$\CI (\Bx_{0}, \Tld{\mu}_{\Bx_{0}},\Tld{\Gl}_{\Bx_{0}})\ge 0$ for $\Tld{\pi}$ almost all 
$\Bx_{0}\in\dOm_{2}$. Thus, we have proved the
inequality \rf{concent-varifold1} for $\Tld{\pi}$ a.e. $\Bx_{0}\in\bra{\GO}$.
This completes the proof of Theorem~\ref{semisuff1}.


\noindent\textbf{Acknowledgments.}
This material is based upon work supported by the National Science Foundation
under Grant No. 0094089. The authors are indebted to Lev Truskinovsky for
sharing his ideas and insights and for many invaluable comments and suggestions.

\begin{thebibliography}{10}

\bibitem{ball89}
John~M. Ball.
\newblock A version of the fundamental theorem for {Y}oung measures.
\newblock In {\em PDEs and continuum models of phase transitions (Nice, 1988)},
  volume 344 of {\em Lecture Notes in Phys.}, pages 207--215. Springer,
  Berlin-New York, 1989.

\bibitem{ball02}
John~M. Ball.
\newblock Some open problems in elasticity.
\newblock In {\em Geometry, mechanics, and dynamics}, pages 3--59. Springer,
  New York, 2002.

\bibitem{bama84}
John~M. Ball and J.~E. Marsden.
\newblock Quasiconvexity at the boundary, positivity of the second variation
  and elastic stability.
\newblock {\em Arch. Rational Mech. Anal.}, 86(3):251--277, 1984.

\bibitem{bamu89}
John~M. Ball and F.~Murat.
\newblock Remarks on {C}hacon's biting lemma.
\newblock {\em Proc. Amer. Math. Soc.}, 107(3):655--663, 1989.

\bibitem{cara29}
C.~Carath{\'{e}}odory.
\newblock \"{U}ber die {V}ariationsrechnung bei mehrfachen {I}ntegralen.
\newblock {\em Acta Math. Szeged}, 4:401--426, 1929.

\bibitem{deDon35}
T.~De~Donder.
\newblock {\em Th{\'e}orie invariantive du clacul des variations}.
\newblock Hayez, Brussels, 1935.

\bibitem{evga92}
Lawrence~C. Evans and Ronald~F. Gariepy.
\newblock {\em Measure theory and fine properties of functions}.
\newblock Studies in Advanced Mathematics. CRC Press, Boca Raton, FL, 1992.

\bibitem{fedebk69}
Herbert Federer.
\newblock {\em Geometric measure theory}.
\newblock Die Grundlehren der mathematischen Wissenschaften, Band 153.
  Springer-Verlag New York Inc., New York, 1969.

\bibitem{fons92}
Irene Fonseca.
\newblock Lower semicontinuity of surface energies.
\newblock {\em Proc. Roy. Soc. Edinburgh Sect. A}, 120(1-2):99--115, 1992.

\bibitem{fomupe98}
Irene Fonseca, Stefan M{\"u}ller, and Pablo Pedregal.
\newblock Analysis of concentration and oscillation effects generated by
  gradients.
\newblock {\em SIAM J. Math. Anal.}, 29(3):736--756 (electronic), 1998.

\bibitem{gitrbk01}
David Gilbarg and Neil~S. Trudinger.
\newblock {\em Elliptic partial differential equations of second order}.
\newblock Classics in Mathematics. Springer-Verlag, Berlin, 2001.
\newblock Reprint of the 1998 edition.

\bibitem{grtr08}
Y.~Grabovsky and L.~Truskinovsky.
\newblock Metastability in nonlinear elsticity.
\newblock {\em To be submitted}.

\bibitem{hest48}
Magnus~R. Hestenes.
\newblock Sufficient conditions for multiple integral problems in the calculus
  of variations.
\newblock {\em Amer. J. Math.}, 70:239--276, 1948.

\bibitem{ks}
R.~V. Kohn and G.~Strang.
\newblock Optimal design and relaxation of variational problems.
\newblock {\em Comm. Pure Appl. Math.}, 39:113--137, 139--182 and 353--377,
  1986.

\bibitem{kris94}
J.~Kristensen.
\newblock Finite functionals and young measures generated by gradients of
  sobolev functions.
\newblock Technical Report Mat-Report No. 1994-34, Mathematical Institute,
  Technical University of Denmark, 1994.

\bibitem{lepa41}
Th. Lepage.
\newblock Sur les champs g\'eod\'esiques des int\'egrales multiples.
\newblock {\em Acad. Roy. Belgique. Bull. Cl. Sci. (5)}, 27:27--46, 1941.

\bibitem{meye65}
Norman~G. Meyers.
\newblock Quasi-convexity and lower semi-continuity of multiple variational
  integrals of any order.
\newblock {\em Trans. Amer. Math. Soc.}, 119:125--149, 1965.

\bibitem{morr52}
Charles~B. Morrey, Jr.
\newblock Quasi-convexity and the lower semicontinuity of multiple integrals.
\newblock {\em Pacific J. Math.}, 2:25--53, 1952.

\bibitem{pedrbk97}
Pablo Pedregal.
\newblock {\em Parametrized measures and variational principles}.
\newblock Progress in Nonlinear Differential Equations and their Applications,
  30. Birkh\"auser Verlag, Basel, 1997.

\bibitem{stein70}
Elias~M. Stein.
\newblock {\em Singular integrals and differentiability properties of
  functions}.
\newblock Princeton Mathematical Series, No. 30. Princeton University Press,
  Princeton, N.J., 1970.

\bibitem{tahe01}
Ali Taheri.
\newblock Sufficiency theorems for local minimizers of the multiple integrals
  of the calculus of variations.
\newblock {\em Proc. Roy. Soc. Edinburgh Sect. A}, 131(1):155--184, 2001.

\bibitem{tart79a}
L.~Tartar.
\newblock Compensated compactness and applications to partial differential
  equations.
\newblock In {\em Nonlinear analysis and mechanics: Heriot-Watt Symposium, Vol.
  IV}, volume~39 of {\em Res. Notes in Math.}, pages 136--212. Pitman, Boston,
  Mass., 1979.

\bibitem{weyl35}
Hermann Weyl.
\newblock Geodesic fields in the calculus of variations of multiple integrals.
\newblock {\em Ann. of Math.}, 36:607--629, 1935.

\bibitem{young33}
L.~C. Young.
\newblock Approximation by polygons in the calculus of variations.
\newblock {\em Proc. Roy. Soc. London, Ser. A}, 141:325--341, 1933.

\bibitem{young37}
L.~C. Young.
\newblock Generalized curves and the existence of an attained absolute minimum
  in the calculus of variations.
\newblock {\em Comptes Rendue, Soc. Sciences \& Lettres, Warsaw, cl. III},
  30:212--234, 1933.

\end{thebibliography}
\def\cprime{$'$}

\end{document}